\documentclass[final,12pt]{elsarticle}
\usepackage{amsfonts,graphics}
\usepackage{amssymb,amsmath}
\usepackage{euscript}
\usepackage{epsfig}
\usepackage{graphicx}
\usepackage{latexsym} 

\def\h2{ { {\mathcal H}_2} }
\def\hinf{ {{\mathcal H}_{\infty}} }

\def\H2{ { {\mathcal H}_2} }
 
\newfont{\Bb}{msbm10 scaled\magstep0}
\def\IR{\mbox {\Bb R}}
\def\IC{\mbox {\Bb C}}

\newcommand{\bfeta}{\mbox{\boldmath$\eta$} }
\newcommand{\bfxi}{\mbox{\boldmath$\xi$} }

\newcommand{\bfSigma}{\mbox{\boldmath$\Sigma$} }
\newcommand{\bfPi}{\mbox{\boldmath$\Pi$} }
\newcommand{\bfGamma}{\mbox{\boldmath$\Gamma$} }
\newcommand{\bfXi}{\mbox{\boldmath$\Xi$} }
\newcommand{\bfPhi}{\mbox{\boldmath$\Phi$} }
\newcommand{\newIRKA}{\textsf{\small IRKA}} 
 \newcommand{\newinIRKA}{\textsf{\small InxIRKA}}

\newcommand{\IL}{{\mathbb{L}}}
\newcommand{\IM}{{\mathbb{M}}}

\newcommand{\bfA}{ {\mathbf A} }
\newcommand{\bfB}{ {\mathbf B} }
\newcommand{\bfC}{ {\mathbf C} }
\newcommand{\bfD}{ {\mathbf D} }
\newcommand{\bfE}{ {\mathbf E} }
\newcommand{\bfF}{ {\mathbf F} }
\newcommand{\bfG}{ {\mathbf G} }
\newcommand{\bfH}{ {\mathbf H} }
\newcommand{\bfI}{ {\mathbf I} }

\newcommand{\bfK}{ {\mathbf K} }
\newcommand{\bfL}{ {\mathbf L} }
\newcommand{\bfM}{ {\mathbf M} }

\newcommand{\bfQ}{ {\mathbf Q} }

\newcommand{\bfV}{ {\mathbf V} }
\newcommand{\bfW}{ {\mathbf W} }
\newcommand{\bfX}{ {\mathbf X} }
\newcommand{\bfY}{ {\mathbf Y} }

\newcommand{\tildebfVr}{ \widetilde{\bfV}_r }
\newcommand{\tildebfWr}{ \widetilde{\bfW}_r }
\newcommand{\tildebfv}{ \widetilde{\bfv}}
\newcommand{\tildebfw}{ \widetilde{\bfw}}

\newcommand{\relnorm}[2]{
\mbox{ {\small $#1 \times 10^{-#2}$ } }  }

\newcommand{\cbfG}{\mbox{\boldmath${\EuScript{G}}$} }

\newcommand{\cbfH}{\mbox{\boldmath${\EuScript{H}}$} }

\newcommand{\hatcbfH}{\mbox{\boldmath$\widehat{\EuScript{H}}$} }
\newcommand{\tildecbfH}{\mbox{\boldmath$\widetilde{\EuScript{H}}$} }
\newcommand{\cbfK}{\mbox{\boldmath${\EuScript{K}}$} }
\newcommand{\hatcbfK}{\mbox{\boldmath$\widehat{\EuScript{K}}$} }
\newcommand{\tildecbfK}{\mbox{\boldmath$\widetilde{\EuScript{K}}$} }

\newcommand{\cbfB}{\mbox{\boldmath${\EuScript{B}}$}}
\newcommand{\cbfM}{\mbox{\boldmath${\EuScript{M}}$}}
\newcommand{\cbfE}{\mbox{\boldmath${\EuScript{E}}$}}
\newcommand{\hatcbfB}{\mbox{\boldmath$\widehat{\EuScript{B}}$}}
\newcommand{\tildecbfB}{\mbox{\boldmath$\widetilde{\EuScript{B}}$}}
\newcommand{\cbfC}{\mbox{\boldmath${\EuScript{C}}$}}
\newcommand{\subcbfC}{\mbox{\scriptsize{\boldmath${\EuScript{C}}$}}}
\newcommand{\subcbfB}{\mbox{\scriptsize{\boldmath${\EuScript{B}}$}}}

\newcommand{\hatcbfC}{\mbox{\boldmath$\widehat{\EuScript{C}}$}}
\newcommand{\tildecbfC}{\mbox{\boldmath$\widetilde{\EuScript{C}}$}}
\newcommand{\cbfP}{\mbox{\boldmath${\EuScript{P}}$}}
\newcommand{\cbfQ}{\mbox{\boldmath${\EuScript{Q}}$}}

\newcommand{\hatcbfQ}{\mbox{\boldmath$\widehat{\EuScript{Q}}$} }
\newcommand{\hatcbfP}{\mbox{\boldmath$\widehat{\EuScript{P}}$} }

\newcommand{\bfb}{ {\mathbf b} }
\newcommand{\bfc}{ {\mathbf c} }

\newcommand{\bfe}{ {\mathbf e} }

\newcommand{\bfu}{ {\mathbf u} }
\newcommand{\bfv}{ {\mathbf v} }
\newcommand{\bfw}{ {\mathbf w} }
\newcommand{\bfx}{ {\mathbf x} }
\newcommand{\bfy}{ {\mathbf y} }

\newcommand{\ignore}[1]{}

\newtheorem{theorem}{Theorem}[section]

\newtheorem{corollary}{Corollary}[section]

\newtheorem{algorithm}{Algorithm}

\begin{document}
\begin{frontmatter}

\title{Inexact Solves in Interpolatory Model Reduction\tnoteref{NSF} \tnoteref{Dedication}}

\author[Bea]{Christopher Beattie}\ead{beattie@vt.edu}
\author[Bea]{Serkan Gugercin}\ead{gugercin@math.vt.edu}    
\author[Bea]{Sarah Wyatt}\ead{sawyatt@vt.edu}               

\address[Bea]{Department of Mathematics,
Virginia Tech, Blacksburg, VA, 24061-0123}  

\tnotetext[Dedication]{
Dedicated to Danny Sorensen on the occasion of his 65th birthday.
}

\tnotetext[NSF]{
This work was supported in part by the NSF through Grants   
DMS-0505971 and DMS- 0645347
}


\begin{abstract}                          
We investigate the  use of inexact solves for
interpolatory model reduction and consider associated
perturbation effects on the underlying model reduction problem. 
We give bounds on system perturbations induced by inexact solves and
relate this to termination criteria for iterative solution methods.  
We show that when a Petrov-Galerkin framework is employed for the inexact solves, 
the associated reduced order model is an exact interpolatory model for a nearby full-order system; 
thus demonstrating backward stability.
We also give evidence that for $\h2$-optimal interpolation points,
interpolatory model reduction is robust with respect to 
perturbations due to inexact solves. 
Finally, we demonstrate the effecitveness of direct use of inexact solves 
in optimal ${\mathcal H}_2$
approximation. The result is an effective model reduction
strategy that is applicable in realistically large-scale settings.
\end{abstract}

\end{frontmatter}

\section{Introduction}
The simulation of dynamical systems constitutes a basic framework for the 
modeling and control
of many complex phenomena of interest in science and
industry. The need for ever greater
model fidelity 
often leads to
computational tasks that make unmanageably large demands on resources.  
Efficient model utilization becomes a critical
consideration in such large-scale problem settings and
motivates the development of strategies for model reduction.

We consider here 
 linear time invariant multi-input/multi-output (MIMO) systems that have 
a state space form (in the Laplace transform domain) as
\begin{align}
\mbox{ Find} \  \widehat{\mathbf{v}}(s)
\  \mbox{such that} \ 
 \cbfK(s)\,\widehat{\mathbf{v}}(s)= \cbfB(s)\widehat{\mathbf{u}}(s) \label{generalMIMO} ,~~
 \  \mbox{then} \  \widehat{\bfy}(s)\stackrel{\mbox{\tiny{def}}}{=}\cbfC(s)\,\widehat{\mathbf{v}}(s). 
\end{align} 
Here, $\widehat{\mathbf{u}}(s)$ and $ \widehat{\bfy}(s)$  denote Laplace-transformed system inputs and outputs, respectively;  $\widehat{\mathbf{v}}(s)$ represents the internal system state.
We assume that $\cbfC(s)\in \IC^{p \times n}$ 
and $\cbfB(s)\in \IC^{n\times m}$ 
are analytic in the right half plane; and that 
$\cbfK(s)\in \IC^{n \times n}$ is analytic and full rank 
throughout the right half plane.  Solving for $\widehat{\bfy}(s)$ in terms of $\widehat{\bfu}(s)$, we obtain
\begin{equation}    \label{HdecompD}
   \widehat{\bfy}(s) = \cbfC(s)\cbfK(s)^{-1}\cbfB(s)   \widehat{\bfu}(s)= \cbfH(s) \widehat{\bfu}(s).
\end{equation}
This representation of the \emph{transfer function}, 
\begin{equation} \label{gencoprime}
\cbfH(s)=\cbfC(s)\cbfK(s)^{-1}\cbfB(s), 
\end{equation} 
we refer to as a \emph{generalized coprime realization}. 
Standard first-order descriptor system realizations, with $ \cbfH(s) = \bfC \left( s \bfE-\bfA\right)^{-1}\bfB$ for constant matrices $\bfE,\ \bfA \in \IR^{n \times n}$,  $\bfB \in \IR^{n\times m}$, and 
$\bfC \in \IR^{p \times n}$ evidently fit this pattern with  $\cbfC(s)= \bfC$, $\cbfB(s)=\bfB$, and $\cbfK(s)=s \bfE-\bfA$.
However, many dynamical systems can be described more naturally with generalized coprime realizations.  For example, a system that includes internal system delays as well as  transmission/propagation delays in its input and output could be described with a model
\begin{equation} \label{ExDelaySS}
\bfE \dot{\bfx}(t)  = \bfA_{0}\, \bfx(t) + \bfA_{1}\, \bfx(t-\tau_{sys}) 
+ \bfB\, \bfu(t-\tau_{\imath np}),  
 \quad \bfy(t)  =  \bfC\, \bfx(t-\tau_{out}) 
\end{equation}
for $\tau_{sys},\,\tau_{\imath np},\,\tau_{out}>0$, and $\bfE,\,\bfA_{0},\, \bfA_{1}  \in \IR^{n \times n}$, $\bfB \in \IR^{n \times m} $
and $\bfC \in \IR^{p \times n}$. Taking the Laplace transformation of (\ref{ExDelaySS}) yields the transfer function 
{\small
$$
\cbfH(s) =  \cbfC(s) \cbfK(s)^{-1} \cbfB(s) =  \left(e^{-s\,\tau_{out}}\,\bfC\right)\left(s\,\bfE -\bfA_{0} -e^{-s\,\tau_{sys}}\, \bfA_{1}\right)^{-1} \left(e^{-s\,\tau_{\imath np}}\,\bfB\right),
$$
}
which has the form of (\ref{gencoprime}).  
  The form of (\ref{gencoprime}) can accomodate greater generality than this, of course, including memory convolution involving higher derivatives, second and higher-order polynomial differential equations, 
systems described via integro-differential equations, and systems where state variables may be coupled through infinite dimensional subsystems (possibly modeling internal propagation or diffusion).  See Table \ref{GCPexampleTable} for other examples and  \cite{beattie2009ipm} for further discussion. 

\begin{table}[htdp]
\caption{Examples of Generalized Coprime System Realizations}
\begin{center}
\begin{tabular}{|c|c|}
\hline
Descriptor Systems & $\bfC(s\bfE-\bfA)^{-1}\bfB$  ($\bfE$ possibly singular)\\
Delay Systems & $(e^{-s\,\tau_{out}}\bfC)(s\bfI-\bfA_0-e^{-s\,\tau_{sys}}\bfA_1)^{-1}(e^{-s\,\tau_{\imath np}}\bfB)$\\
Second Order Systems & $(s\bfC_1+\bfC_0)(s^2\bfM+s\bfG+\bfK)^{-1}\bfB$ \\
Weighted Systems & $\bfW_o(s)\bfC(s\bfI-\bfA)^{-1}\bfB\bfW_\imath(s)$\\
\hline
\end{tabular}
\end{center}
\label{GCPexampleTable}
\end{table}%

In many applications, the state space dimension, $n$, is
too large for efficient system simulation and control computation, so the cases of interest for us
 here have state space dimension 
vastly larger than input and output dimensions: 
$n \gg m, p$.  See \cite{korvink2005obc} for a recent
collection of such benchmark problems.

 The goal is to produce a reduced system
that will have approximately the same response (output) as
the original system for any given input $\bfu(t)$.
For a given reduced-order $r \ll n$,
we construct reduced order models through 
a Petrov-Galerkin approximation of (\ref{generalMIMO}):
Select  full rank  matrices $\bfV_{r}
\in {\mathbb R}^{n \times r}$ and $\bfW_{r}\in {\mathbb R}^{n \times r}$.     For any input, $\bfu(t)$, 
the reduced system output, $\bfy_{r}(t)$, is then defined (in the  Laplace transform domain) as:
\begin{align}
\mbox{ Find} \  \widehat{\mathbf{v}}(s)\in\mbox{Ran}(\bfV_{r}) 
\  \mbox{such that} \  &
\bfW_{r}^T  \left(  \cbfK(s)\,\widehat{\mathbf{v}}(s)- \cbfB(s)\widehat{\mathbf{u}}(s) \right) =0   \\
 \  \mbox{then} \  \widehat{\bfy}_{r}(s)\,\stackrel{\mbox{\tiny{def}}}{=}\,        &
                                \cbfC(s)\,\widehat{\mathbf{v}}(s)\,
\end{align} 
which defines the reduced transfer function as,
\begin{equation}    \label{HredDecompD}
    \cbfH_r(s)=\cbfC_r(s)\cbfK_r(s)^{-1}\cbfB_r(s),
\end{equation}
where
\begin{align}  \label{red_projection} \nonumber
\cbfK_r(s)= \bfW_r^{T} \cbfK(s) \bfV_r  & \in \IC^{r \times r},\quad
\cbfB_r(s) = \bfW_r^{T}\cbfB(s)\in \IC^{r\times m}, \\
& {\rm and}\quad\cbfC_r(s) = \cbfC(s) \bfV_r\in \IC^{p \times r}.
\end{align}

\section{Interpolatory Model Reduction}  Interpolatory reduced order models are designed to exactly reproduce certain system response components that result from
 inputs having specified frequency content and growth.  The approach has been described for standard first-order system realizations in \cite{grimme1997kpm,antoulas2005approximation,gallivan2003mrv,antoulas2010imr} and extended to generalized coprime realizations in  \cite{beattie2009ipm}.
 We summarize the basic elements of this approach below. 
 
A set of points  $\left\{\mu_i\right\}_{i=1}^r\subset \IC $ and (nontrivial) direction vectors 
$\left\{\bfc_i\right\}_{i=1}^r\subset \IC^{p}$ 
constitute left tangential interpolation data for the reduced model,  $\cbfH_r(s)$, if 
\begin{equation}\label{lefttang}
\bfc_i^T \cbfH(\mu_i)=\bfc_i^T \cbfH_r(\mu_i)\quad\mbox{ for each }i=1,\,\ldots,\,r.
\end{equation}
  Likewise,  $\left\{\sigma_j\right\}_{j=1}^r$, and associated directions $\left\{\bfb_j\right\}_{j=1}^r\subset \IC^{m}$, constitute right tangential interpolation data
 for the reduced model,  $\cbfH_r(s)$,  if 
\begin{equation}\label{righttang}
 \cbfH(\sigma_j) \bfb_j = \cbfH_r(\sigma_j) \bfb_j \quad \mbox{ for each }j=1,\,\ldots,\,r.
 \end{equation}
Given left and right tangential interpolating data, interpolatory model
reduction may be implemented by first solving the linear systems: 
\begin{equation}\label{lefttangkrylov}
\mbox{Find }\bfw_i \mbox{ such that } \bfw_i^T\cbfK(\mu_i) =\bfc_i^T\cbfC(\mu_i)\quad \mbox{ for }i=1,\,\ldots,\,r, \mbox{ and}
 \end{equation} 
 \begin{equation}\label{righttangkrylov}
\mbox{find }\bfv_i \mbox{ such that } \cbfK(\sigma_j) \bfv_j = \cbfB(\sigma_j)\bfb_j\quad\mbox{ for  }j=1,\,\ldots,\,r.
 \end{equation}
We assume that the two point sets $\left\{\mu_i\right\}_{i=1}^r$ and $\left\{\sigma_j\right\}_{j=1}^r$ each consist of 
$r$ distinct points and that the vectors $\{\bfv_1,\,\cdots,\,\bfv_r\}$ and 
$\{\bfw_1,\,\cdots,\,\bfw_r\}$ are linearly independent sets.   These vectors constitute ``primitive bases" for the subspaces
$\mathcal{V}_r=\mbox{span}\{\bfv_1,\,\cdots,\,\bfv_r\}$ and 
$\mathcal{W}_r=\mbox{span}\{\bfw_1,\,\cdots,\,\bfw_r\}$. 
Define the associated matrices:
\begin{align}
\hspace*{-1ex}\bfV_r =\  \left[~\bfv_1,~\cdots,~\bfv_r~\right] \  
=&\  \left[~\cbfK(\sigma_1)^{-1}\cbfB(\sigma_1)\bfb_1,\,\cdots,
\,\cbfK(\sigma_r)^{-1}\cbfB(\sigma_r)\bfb_r~\right], \label{k1}\\
\bfW_r^T =
\ \left[\begin{array}{c}
\bfw_1^T\\
\vdots \\
\bfw_r^T
\end{array}\right] = &
\ \left[\begin{array}{c}
\bfc_1^T\cbfC(\mu_1)\cbfK(\mu_1)^{-1}\\
\vdots \\
\bfc_r^T\cbfC(\mu_r)\cbfK(\mu_r)^{-1}
\end{array}\right].  \label{k2}
\end{align}
The reduced model, $\cbfH_r(s)$, as
defined in (\ref{HredDecompD}) and  (\ref{red_projection}) using $\bfV_r $ and $\bfW_r$ 
from (\ref{k1}) and (\ref{k2}),  interpolates $\cbfH(s)$ at the $2r$ points  
$\left\{\mu_i\right\}_{i=1}^r$ and $\left\{\sigma_j\right\}_{j=1}^r$, in respective output 
directions  $\left\{\bfc_i\right\}_{i=1}^r$ and input directions $\left\{\bfb_j\right\}_{j=1}^r$; that is,
conditions (\ref{lefttang}) and (\ref{righttang}) are satisfied.
If 
$\mu_k=\sigma_k$ for some $k$ then first order bitangential moments match as well: 
$$
\bfc_k^T\,\cbfH'(\mu_k)\,\bfb_k = \bfc_k^T\,\cbfH_r'(\mu_k)\,\bfb_k
$$
 Interpolation of higher order derivatives of $\cbfH(s)$ can be accomplished with similar constructions as well;
  see \cite{beattie2009ipm,antoulas2010imr} and references therein.

For large-scale settings with millions of degrees of freedom, 
interpolatory model reduction  has become 
the method of choice since it does 
not require dense matrix operations; the major computational cost lies in solving the (often sparse)
linear systems in (\ref{lefttangkrylov}) and (\ref{righttangkrylov}).
 This contrasts with Gramian-based model reduction approaches such as 
 balanced truncation \cite{mullis1976som,moore1981pca}, optimal 
 Hankel norm approximation \cite{glover1984aoh} and singular perturbation approximation
\cite{liu1989spa} where large-scale Lyapunov equations need to be solved. 
 Moreover, these computational advantages have been enhanced for standard first order state-space realizations by strategies for optimal selection of  tangential interpolation data, see \cite{gugercin2008hmr}.

\subsection{Inexact Interpolatory Model Reduction} \label{inexactkrylov}
The basic framework for interpolatory model reduction presumes
 that the key equations (\ref{lefttangkrylov}) and (\ref{righttangkrylov}) may be solved exactly or nearly so, at least to an accuracy 
 associated with machine precision. 
Direct solution methods, employing sparse
factorization strategies, for example, are capable of handling systems of significantly large order. 
However since the need for ever greater modeling detail and fidelity can drive system order to
the order of millions, the use of direct solvers for the linear
systems  (\ref{lefttangkrylov}) and (\ref{righttangkrylov}) often becomes infeasible and iterative methods must be employed that 
terminate with possibly coarse approximate solutions to the linear systems.  
We consider and evaluate issues related to these approaches here.

Suppose $\{\widehat{\bfv}_1,\,\cdots,\, \widehat{\bfv}_r\}$ and 
$\{\widehat{\bfw}_1,\,\cdots,\, \widehat{\bfw}_r\}$ are linearly 
independent sets in $\IC^n$ and define
\begin{equation}   \label{defineVhatWhat}
\widehat{\bfV}_r =\  \left[~ \widehat{\bfv}_1,~\cdots,~ \widehat{\bfv}_r~\right]
\qquad
\widehat{\bfW}_r^T =
\ \left[\begin{array}{c}
\widehat{\bfw}_1^T\\
\vdots \\
\widehat{\bfw}_r^T
\end{array}\right].
\end{equation}
 $\widehat{\bfw}_i$  and $\widehat{\bfv}_j$ will be viewed as approximate solutions to the linear systems (\ref{lefttangkrylov}) and (\ref{righttangkrylov}) and accordingly we will refer to them as ``inexact" solutions to (\ref{lefttangkrylov}) and (\ref{righttangkrylov}). 
 Nonetheless, unless otherwise stated, these vectors can be any arbitrarily chosen linearly independent vectors in $\IC^n$.     

Define residuals, $\bfxi_i$ and $\bfeta_j$, corresponding to $\widehat{\bfw}_i$  and $\widehat{\bfv}_j$, as
\begin{equation}\label{resb}
\bfxi_i = \cbfK(\mu_i)^T\, \widehat{\bfw}_i  -\cbfC(\mu_i)^T \, \bfc_i \quad
\mbox{ and } \quad
\bfeta_j =  \cbfK(\sigma_j)  \widehat{\bfv}_j - \cbfB(\sigma_j)\bfb_j \ 
\end{equation}
The deviations from the corresponding exact solutions are then
\begin{equation} \label{deltavj} 
{\delta}\bfw_i  = \widehat{\bfw}_i - \bfw_i = \cbfK(\mu_i)^{-T}\, \bfxi_i  \quad \mbox{and}
\quad {\delta}\bfv_j  = \widehat{\bfv}_j - \bfv_j = \cbfK(\sigma_j)^{-1} \, \bfeta_j.
\end{equation} 
The resulting (inexact) basis matrices destined for use in a reduced order model are
\begin{align}
\widehat{\bfW}_r &= \bfW_r +  \left[\delta\bfw_1,~\cdots,~\delta\bfw_r\right]\label{What}\\
\widehat{\bfV}_r &= \bfV_r +  \left[\delta\bfv_1,~\cdots,~\delta\bfv_r\right].\label{Vhat}
\end{align}

Define reduced order maps associated with these inexact bases:
\begin{equation}  \label{approx_red_proj}
\hatcbfK_r(s)= \widehat{\bfW}_r^{T} \cbfK(s) \widehat{\bfV}_r,\quad
\hatcbfB_r(s) = \widehat{\bfW}_r^{T}\cbfB(s), \quad 
\mbox{and}\quad \hatcbfC_r(s) = \cbfC(s) \widehat{\bfV}_r,
\end{equation}
together with the associated inexact reduced order transfer function 
$$
\hatcbfH_r(s)= \hatcbfC_r(s) \hatcbfK_r(s)^{-1} \hatcbfB_r(s). 
$$ 
Notice that we are free to make any choice for bases for the subspaces, $\widehat{\mathcal V}_r$ 
and $\widehat{\mathcal W}_r$, in defining $\hatcbfH_r(s)$; no change in the definition of (\ref{approx_red_proj}) is necessary. As a practical matter, it is generally prudent to choose well conditioned bases in computation.

\section{Forward Error}
\subsection{Interpolation Error }
Inexactness in the solution of the key linear systems (\ref{lefttangkrylov}) and (\ref{righttangkrylov}) produces
a computed reduced order transfer 
function, $\hatcbfH_r(s)$ that no longer interpolates $\cbfH(s)$; typically, 
the reduced order system response will no longer match any component of the full order system response 
at any of the complex frequencies $\left\{\mu_i\right\}_{i=1}^r$ and $\left\{\sigma_i\right\}_{i=1}^r$ that have been specified.   
How much response error has been introduced at these points ?

The particular realization taken for a transfer function can create innate sensitivities to perturbations associated with that representation.  
 Define perturbed transfer functions,  
$$
\cbfH_{\delta \subcbfB}(s)=\cbfC(s)\cbfK(s)^{-1}(\cbfB(s)+\delta\cbfB)\quad \mbox{and}\quad \cbfH_{\delta \subcbfC}(s)=(\cbfC(s)+\delta\cbfC)\cbfK(s)^{-1}\cbfB(s).
$$ 
In discussing perurbations in system response caused by $\delta \cbfB$ and $\delta \cbfC$ at $s=\sigma$,  it is natural to introduce the following quantities:
\begin{align*} 
\mbox{cond}_{\subcbfB}(\cbfH(\sigma))&=\frac{\|\cbfC(\sigma)\cbfK(\sigma)^{-1}\|\,\|\cbfB(\sigma)\|}{\|\cbfH(\sigma)\|}\\
\mbox{cond}_{\subcbfC}(\cbfH(\sigma))&=\frac{\|\cbfC(\sigma)\|\,\|\cbfK(\sigma)^{-1}\cbfB(\sigma)\|}{\|\cbfH(\sigma)\|}
\end{align*} 
to be \emph{condition numbers of the transfer function response}, by way of analogy to the condition number of algebraic linear systems. (Unless otherwise noted, norms will always refer to the Euclidean $2$-norm for vectors or the naturally induced spectral norm for matrices).   It is straightforward to show that these quantities measure the relative sensitivity of the system with respect to perturbations in $\cbfB$ and $\cbfC$, respectively:
\begin{align*} 
\frac{\|\cbfH_{\delta \subcbfB}(\sigma)-\cbfH(\sigma)\|}{\|\cbfH(\sigma)\|}&\leq \mbox{cond}_{\subcbfB}(\cbfH(\sigma))\frac{\|\delta\cbfB\,\|}
{\|\cbfB(\sigma)\|}
\quad \mbox{and} \\
\frac{\| \cbfH_{\delta \subcbfC}(\sigma)-\cbfH(\sigma)\|}{\|\cbfH(\sigma)\|}&\leq 
\mbox{cond}_{\subcbfC}(\cbfH(\sigma))
\frac{\|\,\delta\cbfC\|}{\|\cbfC(\sigma)\|}.
\end{align*}

 For values of $s$ such that 
$\cbfK_r(s)$ and $\hatcbfK_r(s)$ are nonsingular, 
define the matrix-valued functions,
\begin{align} 
\cbfP_{r}(s)  = \cbfK(s)\bfV_r 
 \cbfK_r(s)^{-1}&\bfW_r ^{T},   \quad
   \cbfQ_r(s)  =\bfV_r \cbfK_r(s)^{-1}\,
 \bfW_r ^{T} \cbfK(s), \nonumber \\
  \hatcbfP_{r}(s)  = \cbfK(s) \widehat{\bfV}_r 
 \hatcbfK_r(s)^{-1}\widehat{\bfW}_r ^{T}, &   \quad\mbox{and}\quad
   \hatcbfQ_r(s)  = \widehat{\bfV}_r \hatcbfK_r(s)^{-1}\,
  \widehat{\bfW}_r ^{T} \cbfK(s) \label{Proj_def}
  \end{align}
Where defined, $\cbfP_{r}(s),$ $\cbfQ_{r}(s),$ $\hatcbfP_{r}(s),$ and $\hatcbfQ_{r}(s)$ are differentiable  (indeed, analytic)  
with respect to $s$, having derivatives that satisfy:
\begin{equation} \label{Proj_Deriv}
     \mbox{and}\ 
\begin{array}{c}
  \hatcbfP_{r}'(s)  = \left(\bfI-\hatcbfP_{r}\right) \cbfK'(s) \cbfK(s)^{-1} \hatcbfP_{r} \\[.1in]
  \qquad \hatcbfQ_r'(s)  = \hatcbfQ_r \cbfK(s)^{-1} \cbfK'(s)  \left(\bfI-\hatcbfQ_r\right) 
   \end{array}
  \end{equation}
with equivalent expressions for $\cbfP_{r}'(s)$ and $\cbfQ_{r}'(s)$.  
We will make a series of observations about properties of $\hatcbfP_{r}(s)$ and $\hatcbfQ_{r}(s)$ 
which will have immediately apparent parallels to properties for  $\cbfP_{r}(s)$ and $\cbfQ_{r}(s)$. 
 
Observe first that $\hatcbfP_{r}^2 = \hatcbfP_{r} $ and $\hatcbfQ_{r}^2 = \hatcbfQ_{r} $ so both $\hatcbfP_{r}(s)$ and $\hatcbfQ_{r}(s)$ are skew projectors.  These projectors are of interest because the pointwise error in the transfer function can be expressed as
{\small
\begin{align*}
 \cbfH(s) -\hatcbfH_r(s) = &
  \cbfC(s) \left(  \cbfK(s)^{-1} -  
 \widehat{\bfV}_r \hatcbfK_r(s)^{-1}  \widehat{\bfW}_r^T    \right)\cbfB(s) \\
  = & \cbfC(s) \cbfK(s)^{-1}    \left(\bfI -  \hatcbfP_{r}(s) \right)\cbfB(s). 
\end{align*}
}
Similarly, 
{\small
$$
 \cbfH(s) -\hatcbfH_r(s) 
  = \cbfC(s)  \left(\bfI -  \hatcbfQ_{r}(s) \right) \cbfK(s)^{-1}  \cbfB(s) 
$$
}
and 
{\small
$$
 \cbfH(s) -\hatcbfH_r(s) 
 = \cbfC(s)  \left(\bfI -  \hatcbfQ_{r}(s) \right) \cbfK(s)^{-1}
  \left(\bfI -  \hatcbfP_{r}(s) \right)\cbfB(s). $$
}
The derivative of this last expression can be computed with the aid of (\ref{Proj_Deriv}) and observing 
$\cbfK(s)^{-1} \hatcbfP_{r}(s)=\hatcbfQ_{r}(s) \cbfK(s)^{-1} $:
{\small
\begin{align} \label{error_Deriv}
 \cbfH'(s) -\hatcbfH_r'(s) 
 = & \frac{d}{ds}\left[\cbfC(s)\cbfK(s)^{-1}\right]  \left(\bfI -  \hatcbfP_{r}(s) \right)\cbfB(s) \\
 & \quad +  \cbfC(s) \left(\bfI -  \hatcbfQ_{r}(s) \right) \frac{d}{ds}\left[\cbfK(s)^{-1}\cbfB(s)\right]  \nonumber \\
& \qquad - \cbfC(s) \left(\bfI -  \hatcbfQ_{r}(s) \right)  \frac{d}{ds}\left[\cbfK(s)^{-1}\right]
  \left(\bfI -  \hatcbfP_{r}(s) \right)\cbfB(s).  \nonumber
\end{align}
}

We introduce the following ($s$-dependent) subspaces:
\begin{align*}
{\mathfrak{P}}_r(s)  = \mbox{Ran } \cbfP_{r}(s) = \mbox{Ran }\cbfK(s){\bfV}_r,   \quad &
{\mathfrak{Q}}_r(s) =  \mbox{Ker}\left({\bfW}_r^{T}\cbfK(s)\right)^\perp,  \\
\widehat{\mathfrak{P}}_r(s)  = \mbox{Ran } \hatcbfP_{r}(s) = \mbox{Ran }\cbfK(s) \widehat{\bfV}_r,  & \quad 
\widehat{\mathfrak{Q}}_r(s) =  \mbox{Ker}\left(\widehat{\bfW}_r^{T}\cbfK(s)\right)^\perp,  \\
\mathfrak{B}_m(s)  =  \mbox{Ran }\cbfK(s)^{-1}\cbfB(s), & \qquad 
\mathfrak{C}_p(s)  =  \mbox{Ker}\left(\cbfC(s) \cbfK(s)^{-1}\right)^\perp.
\end{align*}
 $\hatcbfP_{r}(s)$ maps vectors in $\IC^n$ onto 
$\widehat{\mathfrak{P}}_r(s) $
along $\widehat{\mathcal{W}}_r^{\perp}$ and $\hatcbfQ_{r}$  maps vectors in $\IC^n$ onto 
$\widehat{\mathcal{V}}_r $
along $\widehat{\mathfrak{Q}}_r(s)^{\perp}$.

Given two  subspaces of $\IC^n$, say ${\mathcal  M}$ and ${\mathcal N}$, we express the proximity 
of one to the other in terms of the angle between the subspaces, $\Theta({\mathcal M},{\mathcal N})\in [0,\frac{\pi}{2}]$  
 defined  as
      $$
	\sup_{ \bfx\in{\mathcal  M}}\inf_{\,\bfy\in{\mathcal  N}} 
	   \frac{\| \bfy -  \bfx \|}{\| \bfx\|}= \sin\Theta({\mathcal M},{\mathcal N}). 
	   $$ 
$\Theta({\mathcal M},{\mathcal  N})$ is the largest canonical angle 
between ${\mathcal  M}$ 
and a ``closest'' subspace $\widehat {\mathcal  N}$ of 
${\mathcal  N}$ having dimension equal to 
$\dim{\mathcal M}$. Notice that if $\dim {\mathcal  N}< \dim {\mathcal  M}$ then 
$\Theta({\mathcal M},{\mathcal  N})=\frac{\pi}{2}$ and
$\Theta({\mathcal M},{\mathcal  N})=0$ if and only if 
${\mathcal  M}\subset {\mathcal  N}$.   
$\Theta({\mathcal M},{\mathcal  N})$ is asymmetrically defined with respect to ${\mathcal  M}$ and ${\mathcal N}$, however if $\dim {\mathcal  N} = \dim {\mathcal  M}$ 
then 
$\Theta({\mathcal M},{\mathcal  N})=\Theta({\mathcal N},{\mathcal  M})$.  
If $\bfPi_{\mathcal M}$ and $\bfPi_{\mathcal N}$ denote orthogonal projectors onto ${\mathcal M}$ and ${\mathcal N}$, respectively, then $\sin\Theta({\mathcal M},{\mathcal  N})=\|(\bfI -\bfPi_{\mathcal M})\bfPi_{\mathcal N} \|$.

The spectral norm of a skew projector can be expressed in terms of the angle between its range and cokernel \cite{szyld2006tmp}. In particular,
\begin{align}
\| \hatcbfP_{r}(s)\| & = \|\bfI -  \hatcbfP_{r}(s)\| 
=\frac{1}{\cos\Theta(\widehat{\mathfrak{P}}_{r}(s),\, \widehat{\mathcal{W}}_r)} \label{P_normbnd}\\
\| \hatcbfQ_{r}(s)\| & = \|\bfI -  \hatcbfQ_{r}(s)\| 
=\frac{1}{\cos\Theta(\widehat{\mathfrak{Q}}_{r}(s),\, \widehat{\mathcal{V}}_r)}  \label{Q_normbnd}
\end{align}

\begin{theorem} \label{lemma_abs}
Given the full-order model $\cbfH(s) = \cbfC(s)\cbfK(s)^{-1}\cbfB(s)$, 
 interpolation points $\{ \sigma_j\} \subset \IC$, $\{ \mu_i\} \subset \IC$ and 
 corresponding tangential directions, 
$\{\bfb_j\} \subset \IC^m$ and $\{\bfc_i\} \subset \IC^p$,
let the inexact interpolatory reduced model $\hatcbfH_r(s) = \hatcbfC_r(s) \hatcbfK_r(s)^{-1} \hatcbfB_r(s)$ be constructed as defined in (\ref{defineVhatWhat})-(\ref{approx_red_proj}).
The (tangential) interpolation error at $\mu_i$ and $\sigma_j$  is 
{\small
\begin{align}
\frac{\| \hatcbfH_r(\sigma_j) \bfb_j -\cbfH(\sigma_j) \bfb_j \|}{\|\cbfH(\sigma_j) \bfb_j \|}
  \leq &\mbox{cond}_{\subcbfB}(\cbfH(\sigma_j) \bfb_j) \ 
  \frac{\sin\Theta\!\left(\mathfrak{C}_p(\sigma_j),\, \widehat{\mathcal{W}}_r\right)} 
{\cos\Theta\!\left( \widehat{\mathfrak{P}}_r(\sigma_j),\, \widehat{\mathcal{W}}_r\right)} 
\ \frac{ \| \bfeta_j\|}{ \|\cbfB(\sigma_j) \bfb_j\|}  \label{RightTanInterpBnd} \\
\frac{\| \bfc_i^T \hatcbfH_r(\mu_i) - \bfc_i^T\cbfH(\mu_i) \|}{\|\bfc_i^T\cbfH(\mu_i) \|}
  \leq & \mbox{cond}_{\subcbfC}(\bfc_i ^T\cbfH(\mu_i))   
  \frac{\sin\Theta\left(\mathfrak{B}_m(\mu_i),\, \widehat{\mathcal{V}}_r\right)} 
{\cos\Theta\left(\widehat{\mathfrak{Q}}_r(\mu_i),\,\widehat{\mathcal{V}}_r \right)} 
\  \frac{\| \bfxi_i\|}{\| \bfc_i ^T \cbfC(\mu_i) \|}. \label{LeftTanInterpBnd}
 \end{align}
 }

  If $\mu_i=\sigma_i$ then, 
{  \small
  \begin{equation} \label{crossmom}
   |\bfc_i^T \hatcbfH_r(\mu_i) \bfb_i - \bfc_i^T\cbfH(\mu_i) \bfb_i |
  \leq   \frac{ \|\cbfK(\mu_i)^{-1} \| \  \| \bfeta_i\|\, \| \bfxi_i\| } 
{\max\!\left(\cos\Theta\!\left(\widehat{\mathfrak{P}}_r(\mu_i),\, \widehat{\mathcal{W}}_r\right),
\cos\Theta\left(\widehat{\mathfrak{Q}}_r(\mu_i),\,\widehat{\mathcal{V}}_r \right)\right)}.
     \end{equation}
     and
     \begin{align} 
  |\bfc_i^T\cbfH'(\mu_i) \bfb_i -\bfc_i^T\hatcbfH_r^{\,\prime}(\mu_i) \bfb_i | & \leq 
  M\left( \frac{\left\|\bfeta_i \right\|}{\cos\Theta(\widehat{\mathfrak{P}}_r(\mu_i),\, \widehat{\mathcal{W}}_r)}  
  + \frac{\left\| \bfxi_i\right\|}{\cos\Theta\!\left(\widehat{\mathfrak{Q}}_r(\mu_i),\, \widehat{\mathcal{V}}_r\right)} \right. \nonumber\\
   + & \left. \frac{\left\|\bfeta_i \right\|}{\cos\Theta(\widehat{\mathfrak{P}}_r(\mu_i),\, \widehat{\mathcal{W}}_r)}
   \frac{ \left\|\bfxi_i\right\|}{\cos\Theta\!\left(\widehat{\mathfrak{Q}}_r(\mu_i),\, \widehat{\mathcal{V}}_r\right)} \right)  \label{crossmomderiv}
  \end{align}
     }
with $M=max(\left\| \frac{d}{ds}\!\left.\left[\bfc_i^T\cbfC\,\cbfK^{-1}\right]\right|_{\mu_i}\right\|,
\left\|\frac{d}{ds}\!\left.\left[\cbfK^{-1}\cbfB \bfb_i \right] \right|_{\mu_i} \right\|, 
\left\|\frac{d}{ds}\!\left.\left[\cbfK^{-1}\right]\right|_{\mu_i}\right\| )$.
\end{theorem}

\textsc{Proof:} 
From (\ref{deltavj}),
$\widehat{\bfv}_j= \cbfK(\sigma_j)^{-1}(\cbfB(\sigma_j)\bfb_j + \bfeta_j)$, which implies then that 
 $\cbfK(\sigma_j)\widehat{\bfv}_j= \cbfB(\sigma_j)\bfb_j + \bfeta_j \in 
\widehat{\mathfrak{P}}_r(\sigma_j) $ and 
 $\left(\bfI - \hatcbfP_{r}(\sigma_j)\right)\left(\cbfB(\sigma_j)\bfb_j + \bfeta_j\right)=0$, which may 
 be rearranged to obtain
 \begin{equation} \label{b2eta}
 \left(\bfI - \hatcbfP_{r}(\sigma_j)\right)\cbfB(\sigma_j)\bfb_j 
 =-\left(\bfI - \hatcbfP_{r}(\sigma_j)\right) \bfeta_j.
  \end{equation}
 Let  $\widehat{\bfPi}$ be the
orthogonal projector taking $\IC^n$ onto  
$\widehat{\mathcal{W}}_r =\mbox{Ker}\left(\hatcbfP_{r}(s) \right)^{\perp}.$ 
One may directly verify that
$
\bfI -  \hatcbfP_{r}(s)=\left(\bfI -  \widehat{\bfPi}\right)\left(\bfI -  \hatcbfP_{r}(s) \right),
$
and
 \begin{align} 
  \hatcbfH_r(\sigma_j) \bfb_j -\cbfH(\sigma_j) \bfb_j 
 & = -\cbfC(\sigma_j)\cbfK(\sigma_j)^{-1}  
 \left(\bfI -  \hatcbfP_{r}(\sigma_j) \right)\cbfB(\sigma_j)\bfb_j  \nonumber \\
 & = \cbfC(\sigma_j)\cbfK(\sigma_j)^{-1}  \left(\bfI -  \hatcbfP_{r}(\sigma_j) \right) \bfeta_j 
 \label{RightTangBackError} \\
 & = \cbfC(\sigma_j)\cbfK(\sigma_j)^{-1} 
  \left(\bfI -  \widehat{\bfPi} \right) \left(\bfI -  \hatcbfP_{r}(\sigma_j) \right)
 \bfeta_j. \nonumber
\end{align} 
Now suppose $\bfGamma$ is an orthogonal projector onto $\mathfrak{C}_p(\sigma_j)$.  We have then that 
$\mbox{Ran}(\bfI-\bfGamma)=\mbox{Ker}(\cbfC(\sigma_j)\cbfK(\sigma_j)^{-1}$, so that 
$\cbfC(\sigma_j)\cbfK(\sigma_j)^{-1})=\cbfC(\sigma_j)\cbfK(\sigma_j)^{-1}\bfGamma$ and
$$
\hatcbfH_r(\sigma_j) \bfb_j -\cbfH(\sigma_j) \bfb_j =
 \cbfC(\sigma_j)\cbfK(\sigma_j)^{-1} \bfGamma
  \left(\bfI -  \widehat{\bfPi} \right) \left(\bfI -  \hatcbfP_{r}(\sigma_j) \right)
 \bfeta_j.
$$

Taking norms, we obtain an estimate yielding (\ref{RightTanInterpBnd}):
{\small
\begin{align*} 
\|  \hatcbfH_r(\sigma_j) \bfb_j -\cbfH(\sigma_j) \bfb_j \|
 & \leq \|\left(\bfI -  \widehat{\bfPi} \right)\bfGamma \left(\cbfC(\sigma_j)\cbfK(\sigma_j)^{-1}\right)^T  \| \cdot   \|\bfI -  \hatcbfP_{r}(\sigma_j) \|\cdot
  \|\bfeta_j\|   \\
& \leq \|\cbfC(\sigma_j)\cbfK(\sigma_j)^{-1}\|\cdot
\frac{\sin\Theta\!\left(\mathfrak{C}_p(\sigma_j),\,\widehat{\mathcal{W}}_r\right)} 
{\cos\Theta(\widehat{\mathfrak{P}}_r(\sigma_j),\, \widehat{\mathcal{W}}_r)}\cdot  \| \bfeta_j\|
\end{align*} 
}
(\ref{LeftTanInterpBnd}) is shown similarly, noting first that 
\begin{equation} \label{c2xi}
\bfc_i^T\cbfC(\mu_i) \left(\bfI - \hatcbfQ_{r}(\mu_i)\right)
 =-\bfxi_i^T\left(\bfI - \hatcbfQ_{r}(\mu_i)\right).
\end{equation}
Defining the orthogonal projector, $\widehat{\bfXi}$, that takes  $\IC^n$ onto  
{\small $\widehat{\mathcal{V}}_r =\mbox{Ran}\left(\hatcbfQ_{r}(s) \right),$}  one observes next 
$ \bfI-\hatcbfQ_{r}(s)= \left(\bfI-\hatcbfQ_{r}(s)\right) \left(\bfI-\widehat{\bfXi}\right)$ so that
\begin{align*} 
\|  \bfc_i^T\hatcbfH_r(\mu_i) - \bfc_i^T\cbfH(\mu_i) \|
 & =  \| \bfc_i^T\cbfC(\mu_i)  
 \left(\bfI -  \hatcbfQ_{r}(\mu_i) \right)\cbfK(\mu_i)^{-1}\cbfB(\mu_i) \| \\
 & \leq  \| \bfxi_i^T
 \left(\bfI -  \hatcbfQ_{r}(\mu_i) \right)\left(\bfI -  \widehat{\bfXi} \right)\cbfK(\mu_i)^{-1}\cbfB(\mu_i) \|   \\
 & \leq \| \bfxi_i\|\, \cdot  \, \|\bfI -  \hatcbfQ_{r}(\mu_i) \|\,\cdot\,
  \|\left(\bfI -  \widehat{\bfXi} \right) \cbfK(\mu_i)^{-1} \cbfB(\mu_i) \| \\
  & \leq \|\cbfK(\mu_i)^{-1} \cbfB(\mu_i) \|\cdot
\frac{\sin\Theta\!\left(\mathfrak{B}_m(\mu_i),\,\widehat{\mathcal{V}}_r\right)} 
{\cos\Theta\!\left(\widehat{\mathfrak{Q}}_r(\mu_i),\, \widehat{\mathcal{V}}_r\right)}\cdot  \| \bfxi_i\|
\end{align*} 
When $\mu_i=\sigma_i$, we have
{\small 
\begin{align*} 
\bfc_i^T\cbfH(\mu_i) \bfb_i - \bfc_i^T \hatcbfH_r(\mu_i) \bfb_i 
& =  \bfc_i^T\cbfC(\mu_i) \left(\bfI -  \hatcbfQ_{r}(\mu_i) \right)\cbfK(\mu_i)^{-1}
  \left(\bfI -  \hatcbfP_{r}(\mu_i) \right)\cbfB(\mu_i) \bfb_i \\
& =  \bfxi_i^T  \left(\bfI -  \hatcbfQ_{r}(\mu_i) \right)\cbfK(\mu_i)^{-1}
  \left(\bfI -  \hatcbfP_{r}(\mu_i) \right) \bfeta_i\\
   & \quad = \left\{\begin{array}{lr}
     \bfxi_i^T \cbfK(\mu_i)^{-1}
  \left(\bfI -  \hatcbfP_{r}(\mu_i) \right) \bfeta_i, & \mbox{or} \\
    \bfxi_i^T  \left(\bfI -  \hatcbfQ_{r}(\mu_i) \right)\cbfK(\mu_i)^{-1} \bfeta_i,
    \end{array} \right.
\end{align*} 
}
leading then to two estimates:
$$
|\bfc_i^T\cbfH(\mu_i) \bfb_i - \bfc_i^T \hatcbfH_r(\mu_i) \bfb_i |
 \leq  \|\bfxi_i\|\cdot  \|\bfeta_i \| \cdot 
 \|\cbfK(\mu_i)^{-1}\|\, \cdot\, \| \bfI -  \hatcbfP_{r}(\mu_i)\| 
 $$
 and
 $$
  |\bfc_i^T\cbfH(\mu_i) \bfb_i - \bfc_i^T \hatcbfH_r(\mu_i) \bfb_i |
 \leq  \|\bfxi_i\|\cdot  \|\bfeta_i \| \cdot 
 \|\cbfK(\mu_i)^{-1}\|\, \cdot\, \| \bfI -  \hatcbfQ_{r}(\mu_i)\|.
$$
These can be combined to yield (\ref{crossmom}).

The last inequality comes from using (\ref{error_Deriv})  with $s=\mu_i$:
\begin{align*}
\bfc_i^T \cbfH'(\mu_i) \bfb_i - \bfc_i^T\hatcbfH_r'(\mu_i) \bfb_i & 
 =  \frac{d}{ds}\!\left.\left[\bfc_i^T\cbfC\,\cbfK^{-1}\right]\right|_{\mu_i}  
 \left(\bfI -  \hatcbfP_{r}(\mu_i) \right)\cbfB(\mu_i) \bfb_i \\
  + \bfc_i^T& \cbfC(\mu_i) \left(\bfI -  \hatcbfQ_{r}(\mu_i) \right)  \frac{d}{ds}\!\left.\left[\cbfK^{-1}\cbfB \bfb_i\right]\right|_{\mu_i}  \\
 - \bfc_i^T\cbfC(\mu_i) & \left(\bfI -  \hatcbfQ_{r}(\mu_i) \right)  \frac{d}{ds}\!\left.\left[\cbfK^{-1}\right]\right|_{\mu_i}  \left(\bfI -  \hatcbfP_{r}(\mu_i) \right)\cbfB(\mu_i) \bfb_i. 
\end{align*}
Then from (\ref{b2eta}),  (\ref{c2xi}), and the Cauchy-Schwarz inequality
{\small
\begin{align*}
& \left|\bfc_i^T \cbfH'(\mu_i) \bfb_i - \bfc_i^T\hatcbfH_r'(\mu_i) \bfb_i\right|
  \leq  \left| \frac{d}{ds}\!\left.\left[\bfc_i^T\cbfC\,\cbfK^{-1}\right]\right|_{\mu_i}  
 \left(\bfI -  \hatcbfP_{r}(\mu_i) \right) \bfeta_i \right| \\
& \qquad \quad  + \left| \bfxi_i^T \left(\bfI -  \hatcbfQ_{r}(\mu_i) \right) \frac{d}{ds}\!\left.\left[\cbfK^{-1}\cbfB \bfb_i \right] \right|_{\mu_i} \right| \\
& \qquad \qquad + \left|\bfxi_i^T \left(\bfI -  \hatcbfQ_{r}(\mu_i) \right)  \frac{d}{ds}\!\left.\left[\cbfK^{-1}\right]\right|_{\mu_i}  \left(\bfI -  \hatcbfP_{r}(\mu_i) \right) \bfeta_i \right|\\
&\quad \leq  \left\| \frac{d}{ds}\!\left.\left[\bfc_i^T\cbfC\,\cbfK^{-1}\right]\right|_{\mu_i}\right\|\cdot
 \frac{\left\|\bfeta_i \right\|}{\cos\Theta(\widehat{\mathfrak{P}}_r(\mu_i),\, \widehat{\mathcal{W}}_r)}\\
& \qquad \quad  + \frac{\left\| \bfxi_i\right\|}{\cos\Theta\!\left(\widehat{\mathfrak{Q}}_r(\mu_i),\, \widehat{\mathcal{V}}_r\right)} \cdot \left\|\frac{d}{ds}\!\left.\left[\cbfK^{-1}\cbfB \bfb_i \right] \right|_{\mu_i} \right\| \\
& \qquad \qquad +  \left\|\frac{d}{ds}\!\left.\left[\cbfK^{-1}\right]\right|_{\mu_i}\right\| \cdot \frac{\left\|\bfeta_i \right\|}{\cos\Theta(\widehat{\mathfrak{P}}_r(\mu_i),\, \widehat{\mathcal{W}}_r)}\frac{ \left\|\bfxi_i\right\|}{\cos\Theta\!\left(\widehat{\mathfrak{Q}}_r(\mu_i),\, \widehat{\mathcal{V}}_r\right)}
\end{align*}
}
which yields the conclusion.   $\Box$

Consider the effect of solving (\ref{lefttangkrylov}) and (\ref{righttangkrylov}) approximately with successively increasing levels of accuracy that force  the residual norms to zero, $\|\bfeta_j\|\rightarrow 0$ and $\| \bfxi_i\|\rightarrow 0$.
The multiplicative behavior of the error bound (\ref{crossmom}) with respect to $\|\bfeta_j\|$ and $\| \bfxi_i\|$ contrasts with the additive behavior seen in (\ref{RightTanInterpBnd}) and (\ref{LeftTanInterpBnd}) and suggests some potential benefit in using the same interpolation points for both left and right interpolation, 
i.e., choosing  $\mu_i=\sigma_i$ for $i=1,\,\ldots,\,r$.  Note that this choice also forces convergent (bitangential) derivative interpolation as 
shown in (\ref{crossmomderiv}). 
Indeed,  choosing $\mu_i=\sigma_i$ for $i=1,\,\ldots,\,r$ is a \emph{necessary} condition for forming $\H2$-optimal interpolatory reduced order models for first-order descriptor realizations,  as we discuss in \S \ref{sec:optimalH2} (see also \cite{gugercin2008hmr}).   Beyond this, there can be notable computational advantages
in choosing $\mu_i=\sigma_i$, since the linear systems to be solved in (\ref{lefttangkrylov}) and (\ref{righttangkrylov}) then have the same coefficient matrix; allowing one potentially to reuse factorizations and preconditioners. 

Certain applications require the retention of structural properties such as symmetry in passing from $\cbfK$ to $\hatcbfK_r$ and one is compelled to choose 
 $\widehat{\bfW}_r = \widehat{\bfV}_r$ (``one-sided" model reduction), so the vectors $\{\widehat{\bfw}_1,\,\cdots,\, \widehat{\bfw}_r\}$ might not be 
 approximate solutions to (\ref{lefttangkrylov}) in the usual sense.  Nonetheless,  the behavior  of the interpolation error is still governed by (\ref{RightTanInterpBnd}) and (\ref{LeftTanInterpBnd}).  We explore this in the following numerical example.

We illustrate the character of the results given in Theorem~\ref{lemma_abs}, bounding the response error at the nominal interpolation points caused by inexact solves in (\ref{lefttangkrylov}) and (\ref{righttangkrylov}).  
 To this end, we consider a delay differential equation of the form introduced in (\ref{ExDelaySS}) taking $n=2000$, $m=p=1$ and $\tau_{\imath np} = \tau_{out} = 0$.  The coefficient matrices for the full order model in (\ref{ExDelaySS}) were taken from \cite{beattie2009ipm}.
We construct multiple reduced models all of order $r=3$ , solving  (\ref{lefttangkrylov}) and (\ref{righttangkrylov}) with different levels of accuracy.  We chose three logarithmically spaced values, $\sigma_1=0.001,\,\sigma_2=0.0316,\,\sigma_3=1.0$,  and fixed them as interpolation points.   We then obtained approximate solutions of varying accuracy to (\ref{lefttangkrylov}) and (\ref{righttangkrylov}) in a manner described in more detail below, assembled the inexact interpolation basis
matrices, $\widehat{\bfV}_r$ and $\widehat{\bfW}_r$,
and obtained reduced models of order $r=3$ having the same internal  delay structure as the original system:
\begin{align*}
\hatcbfH_r(s) &=  \hatcbfC_r(s)\hatcbfK_r(s)^{-1} \hatcbfB_r(s)\\
& =  \bfC\widehat{\bfV}_r\left(s\,\widehat{\bfW}_r^T\bfE \widehat{\bfV}_r-\widehat{\bfW}_r^T\bfA_{0}\widehat{\bfV}_r -e^{-s\,\tau_{sys}}\, \widehat{\bfW}_r^T\bfA_{1}\widehat{\bfV}_r\right)^{-1} \widehat{\bfW}_r^T\bfB
\end{align*}
We considered both the usual ``two-sided" model reduction process that involves approximate solution of both (\ref{lefttangkrylov}) and (\ref{righttangkrylov}) and the ``one-sided" process that involves approximate solutions only to (\ref{righttangkrylov}) to generate $\widehat{\bfV}_r$ and then assigning $\widehat{\bfW}_r = \widehat{\bfV}_r$. 
Linear systems were solved with GMRES
 terminating with a final relative residual below a uniform tolerance denoted by $\varepsilon$. 
 
 We generated reduced order models in this way, varying the relative residual tolerance $\varepsilon$ from $10^{-1}$ down to $10^{-8}$.
  Figure \ref{figure:delay_sigma1} below shows the resulting interpolation errors  $|\cbfH(\sigma_1) - \hatcbfH_r(\sigma_1)|$ 
and bounds from equations (\ref{RightTanInterpBnd}) and (\ref{crossmom})
for one-sided  and two-sided cases, respectively, as $\varepsilon$  varies. 
Observe that the bounds in Theorem~\ref{lemma_abs} predict the convergence behavior of the true error quite well;  the rates (slopes) are matched almost exactly. 
Note also that the interpolation error decays much faster for two-sided reduction than for one-sided reduction 
Indeed,  the ratio of the two errors is close to $\varepsilon$, i.e., for a given tolerance $\varepsilon$, the interpolation error for two-sided reduction is  approximately $\varepsilon$ times smaller than the interpolation error for one-sided reduction.
\begin{figure}[ht]
  \centerline{
\includegraphics[width=10.0cm]{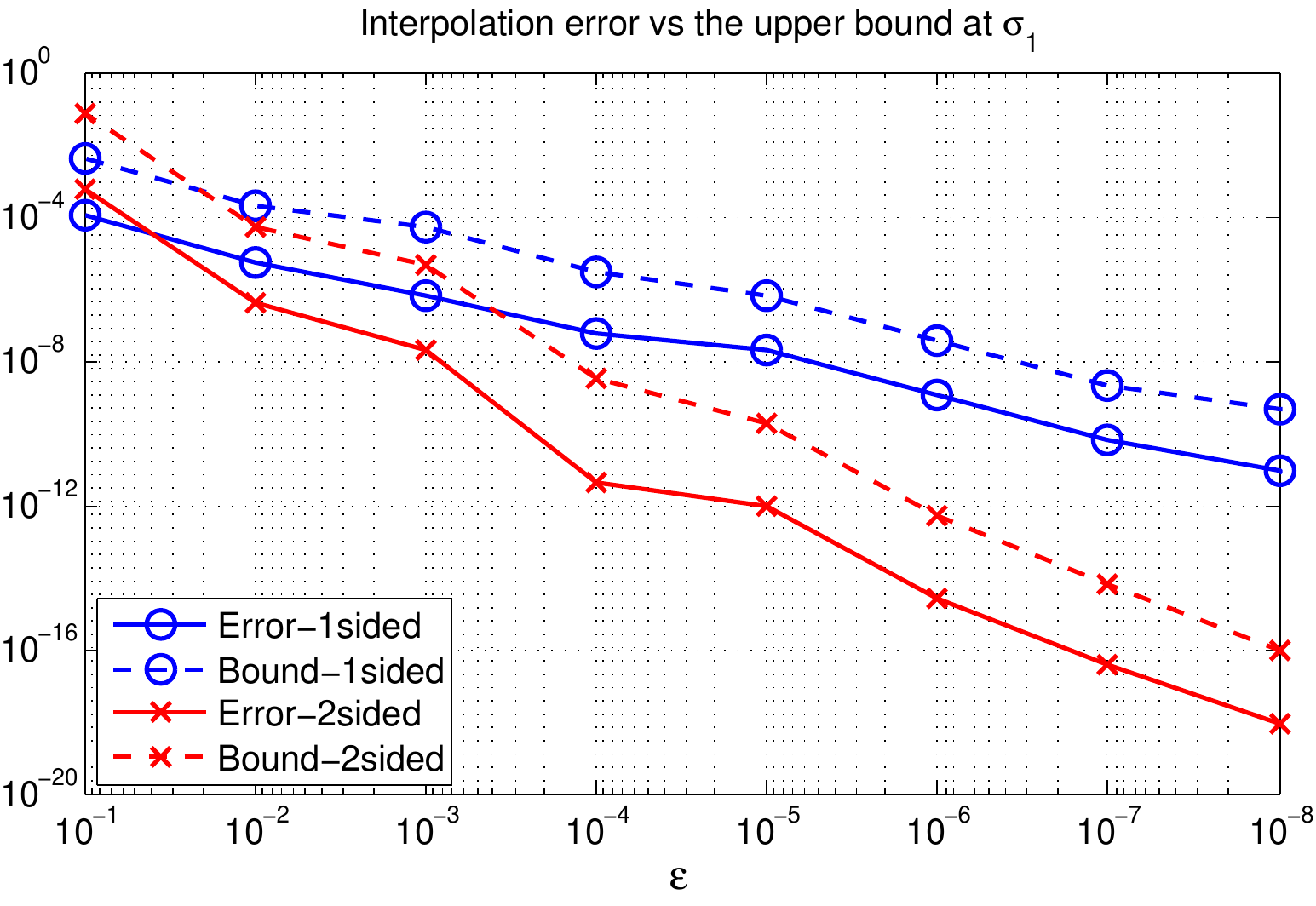}
   }
\caption{Behavior of interpolation error and upper bounds vs $\varepsilon$}  
\label{figure:delay_sigma1}
\end{figure}
  
Analogous results regarding behavior of the bounds and interpolation error are observed
at $\sigma_2$ and $\sigma_3$ and so are omitted for brevity. 


  \subsection{Global Error Bounds} 
  
Thus far we have focussed on the extent to which  interpolation properties are lost in the computed reduced models when inexact solves are introduced into the process, considering in effect \emph{local} error bounds.   
Clearly, it is important to understand the effect of inexact solves on the overall \emph{global} quality of the reduced order model.  There are two commonly used measures for closeness of two conforming dynamical systems (i.e., those with the same input and output dimensions):
\begin{align*}
&\mbox{the $\h2$-norm:} \qquad \|\cbfH-\cbfG\|_{\h2}=\frac{1}{2\pi}\int_{-\infty}^{\infty} \|\cbfH(\imath\omega)-\cbfG(\imath\omega)\|_F^2\,d\omega\\
&\mbox{the $\hinf$-norm: } \qquad \|\cbfH-\cbfG\|_{\hinf}=\max_{\omega\in\IR}\left\|\cbfH(\imath\omega)-\cbfG(\imath\omega)\right\|_2.
\end{align*}

Since  reduced models are completely determined by the subspaces, $\mathcal{V}_r$ and  $\mathcal{W}_r$,
as shown in (\ref{red_projection}),  we first evaluate (in Theorem~\ref{thm:subspacebound}) how much  inexact
interpolatory subspaces, $\widehat{\mathcal{V}}_r$ and $\widehat{\mathcal{W}}_r$, can deviate from  the corresponding true subspaces, ${\mathcal{V}}_r$ and ${\mathcal{W}}_r$, as a result of inexact solves.   The effect of this deviation
on the resulting model reduction (forward) error will be shown in Theorem~\ref{globalFwdBound}.  In this way, we are able to connect model reduction error to observable quantities that are associated with inexact solves, such as the relative stopping criterion $\varepsilon$.

  \begin{theorem} \label{thm:subspacebound}
  Let the columns of $\bfV_r$ and $\widehat{\bfV}_r$ be exact and approximate solutions to (\ref{righttangkrylov}) and the columns of $\bfW_r$ and $\widehat{\bfW}_r$ be exact and approximate solutions to  (\ref{lefttangkrylov}).  Suppose approximate solutions are computed to a relative residual tolerance of $\varepsilon>0$, so that 
  $\|\bfeta_i\| \leq \varepsilon\|\cbfB(\sigma_i)\bfb_i\|$ and $\|\bfxi_i\| \leq \varepsilon\ \|\cbfC(\mu_i)^T \, \bfc_i\|$, where the residuals $\bfeta_i $ and $\bfxi_i$ are defined in (\ref{resb}).
  
Denoting the associated subspaces as $\mathcal{V}_{r}$, $\widehat{\mathcal{V}}_{r}$, $\mathcal{W}_{r}$ and $\widehat{\mathcal{W}}_{r}$ then 
    \begin{align} 
   \sin\Theta(\widehat{\mathcal{V}}_{r},\,\mathcal{V}_r)  & \leq 
  \frac{\varepsilon\ \sqrt{r}}
  {\varsigma_{\min}(\widehat{\bfV}_r\bfD_v)} \label{subspBoundV} \\
  \sin\Theta(\widehat{\mathcal{W}}_{r},\,\mathcal{W}_r)  & \leq 
  \frac{\varepsilon\ \sqrt{r}}
  {\varsigma_{\min}(\widehat{\bfW}_r\bfD_w)} \label{subspBoundW}
       \end{align}
  where $\bfD_v$ and $\bfD_w$ are diagonal scaling matrices defined as
  {\small
 \begin{align*} 
  \bfD_v&=\mbox{diag}\left((\| \cbfK(\sigma_1)^{-1}\|\,\|\cbfB(\sigma_1)\mathbf{b}_1\|)^{-1},\,\ldots,\, (\| \cbfK(\sigma_r)^{-1}\|\,\|\cbfB(\sigma_r)\mathbf{b}_r\|)^{-1}\right)\mbox{ \textrm{and}}\\
   \bfD_w&=\mbox{diag}\left((\| \cbfK(\mu_1)^{-1}\|\,\|\cbfC(\mu_1)^T \, \bfc_1\|)^{-1},\,\ldots,\, (\| \cbfK(\mu_r)^{-1}\|\,\|\cbfC(\mu_r)^T \, \bfc_r\|)^{-1}\right)
  \end{align*} 
  }
and $\varsigma_{\min}(\bfM)$ denotes the smallest singular value of the matrix $\bfM$. 
\end{theorem}
\textsc{Proof:}  We prove (\ref{subspBoundV}).   The proof of (\ref{subspBoundW}) is similar. 

Write $\widehat{\mathbf{V}}_{r} = \mathbf{V}_{r} + \bfE$ with $\bfE = [\bfK(\sigma_1)^{-1} \bfeta_1, \ldots , \bfK(\sigma_r)^{-1} \bfeta_r]$. 
Then 
  {\small
  \begin{align*}
  \sin\Theta(\widehat{\mathcal{V}}_{r},\,\mathcal{V}_r)& = \max_{\hat{\bfv}\in \widehat{\mathcal{V}}_{r}} \min_{\bfv\in \mathcal{V}_r} \frac{\| \bfv-\hat{\bfv}\|}{\| \hat{\bfv} \|} \\
   = \max_{x_i} &\min_{z_i} \frac{\|\sum_{i=1}^r z_i\cbfK(\sigma_i)^{-1}\cbfB(\sigma_i)\mathbf{b}_i - 
  \sum_{i=1}^r x_i \widehat{\bfv}_i\|}{\|\sum_{i=1}^r x_i \widehat{\bfv}_i\|} \\
 =& \max_{x_i} \min_{z_i} \frac{\|\sum_{i=1}^r(z_i-x_i)\cbfK(\sigma_i)^{-1}\cbfB(\sigma_i)\mathbf{b}_i - x_i\cbfK(\sigma_i)^{-1}\bfeta_i\|}
 {\|\sum_{i=1}^r x_i\widehat{\bfv}_i\|} \\
  &  \leq  \max_{x_i}  \frac{\|\sum_{i=1}^r x_i\cbfK(\sigma_i)^{-1}\bfeta_i\|}
 {\|\sum_{i=1}^r x_i\widehat{\bfv}_i\|}  =\max_{\bfx}\frac{\|\bfE\bfx\|}{\| \widehat{\bfV}_r\bfx\|}
 =\max_{\bfx}\frac{\|\bfE\bfD\bfx\|}{\| \widehat{\bfV}_r\bfD\bfx\|}
 \end{align*}
 }
 where $\bfD=\mbox{diag}(d_1,\,\ldots,\,d_r)$ is a diagonal matrix with positive diagonal entries, $d_i>0$, 
 that are fixed but for the moment unspecified. 
 
 Note that 
   \begin{align*}
 \|\bfE\bfD\bfx\|\,\leq&\, \|\bfE\bfD\|\, \|\bfx\|\, \leq\, \sqrt{r}\,\|\bfx\| \,\max_i\left(d_i\| \cbfK(\sigma_i)^{-1}\bfeta_i\|\right) \\
&  \leq \sqrt{r}\, \|\bfx\|\, \max_i\left(d_i\| \cbfK(\sigma_i)^{-1}\|\,\|\bfeta_i\|\right)
   \end{align*}
   Thus we have, 
   {\small
 \begin{equation} \label{SluisBound}
  \sin\Theta(\widehat{\mathcal{V}}_{r},\,\mathcal{V}_r)
  \leq \sqrt{r}\, \frac{\max_i\left(d_i\| \cbfK(\sigma_i)^{-1}\|\,\|\bfeta_i\|\right)}
  {\min_{\bfx}\left(\| \widehat{\bfV}_r\bfD\bfx\|/\|\bfx\|\right)}
  =\sqrt{r}\, \frac{\max_i\left(d_i\| \cbfK(\sigma_i)^{-1}\|\,\|\bfeta_i\|\right)}
  {\varsigma_{\min}(\widehat{\bfV}_r\bfD)}
   \end{equation}
   }
 This bound is valid for any choice of diagonal scalings, $\bfD$, so we can minimize the right hand side of (\ref{SluisBound}) with respect to $d_1,\,\ldots,\,d_r$.  The \emph{Column Equilibration Theorem} of van der Sluis \cite{vanderSluis1969} asserts that the optimal choice of $d_1,\,\ldots,\,d_r$ is such that $d_i\| \cbfK(\sigma_i)^{-1}\|\,\|\bfeta_i\|=C$, independent of $i=1,\,\ldots,\,r$. If inexact solves terminate with residuals satisfying $\|\bfeta_i\|\approx\varepsilon\, \|\cbfB(\sigma_i)\mathbf{b}_i\|$ then we may take
 $C=\varepsilon$ and $d_i=\left(\| \cbfK(\sigma_i)^{-1}\|\,\|\cbfB(\sigma_i)\mathbf{b}_i\|\right)^{-1}$ to achieve the best bound possible with the information given.   This leads to (\ref{subspBoundV}).$\Box$
 
 \vspace{.1in}
As a practical matter, the column scalings used in (\ref{subspBoundV}) and (\ref{subspBoundW}) will not be computationally feasible in realistic settings.  If instead we scale the columns of $\widehat{\bfV}_r$ and $\widehat{\bfW}_r$ to have unit norm (cheap !) --- taking $\widetilde{\bfD}_v= \mbox{diag}\left(1/\| \widehat{\bfv}_1\|,\,\ldots,\, 1/\| \widehat{\bfv}_r\|\right)$ and $\widetilde{\bfD}_w = \mbox{diag}\left(1/\| \widehat{\bfw}_1\|,\,\ldots,\, 1/\| \widehat{\bfw}_r\|\right) $, the bound for (\ref{subspBoundV}) degrades to 
$$
 \sin\Theta(\widehat{\mathcal{V}}_{r},\,\mathcal{V}_r)   \leq \max_i\kappa_2\left(\cbfK(\sigma_i), \widehat{\bfv}_i\right)
  \frac{\varepsilon\ \sqrt{r}}
  {\varsigma_{\min}(\widehat{\bfV}_r \widetilde{\bfD}_v)} 
$$
where $\displaystyle \kappa_2\left(\cbfK(\sigma_i), \widehat{\bfv}_i\right)=\frac{\| \cbfK(\sigma_i)^{-1}\|\,\|\cbfB(\sigma_i)\mathbf{b}_i\|}{\|\widehat{\bfv}_i\|}>1$ is the \emph{condition number of the linear system} (\ref{righttangkrylov}).  A similar expression holds for $ \sin\Theta(\widehat{\mathcal{W}}_{r},\,\mathcal{W}_r) $.  In many cases, these condition numbers have only modest magnitude and the bounds (\ref{subspBoundV}) and (\ref{subspBoundW}) remain descriptive. 

   \begin{theorem} \label{globalFwdBound}
  Let the columns of $\bfV_r$ and $\widehat{\bfV}_r$ be exact and approximate solutions to (\ref{righttangkrylov}) and the columns of $\bfW_r$ and $\widehat{\bfW}_r$ be exact and approximate solutions to  (\ref{lefttangkrylov}).   Let the associated subspaces be denoted as $\mathcal{V}_{r}$, $\widehat{\mathcal{V}}_{r}$, $\mathcal{W}_{r}$ and $\widehat{\mathcal{W}}_{r}$ and the associated reduced order systems be denoted as $\cbfH_r(s)$ (exact) and $\ \hatcbfH_r(s)$ (inexact).  Then 
{\small
 \begin{equation*} 
 \frac{\|\cbfH_r-\hatcbfH_r\|_{\hinf}}{\frac{1}{2}\left( \|\cbfH_{r}\|_{\hinf}+ \|\hatcbfH_{r}\|_{\hinf}\right)}
  \leq\,M\,\max\left(\sin\Theta(\widehat{\mathcal{V}}_{r},\,\mathcal{V}_r),\, \sin\Theta(\widehat{\mathcal{W}}_{r},\,\mathcal{W}_r)\right),
  \end{equation*} 
 }
where 
{\small 
$$ 
M= 2\,\max\left(\frac{\max_{\omega\in\IR}
\mbox{cond}_{\subcbfC}(\cbfH_r(\imath\omega))}
  {\min_{\omega\in\IR}\cos\Theta(\widehat{\mathfrak{Q}}_{r}(\imath \omega),\,\widehat{\mathcal{V}}_r)},
   \frac{\max_{\omega\in\IR} \mbox{cond}_{\subcbfB}(\hatcbfH_r(\imath\omega))}
    {\min_{\omega\in\IR}\cos\Theta(\mathfrak{P}_{r}(\imath\omega),\,\mathcal{W}_r)} 
    \right)
    $$
    } and
    \begin{align*} 
\mbox{cond}_{\subcbfB}(\hatcbfH_r(s))&=\frac{\|\hatcbfC_r(s)\hatcbfK_r(s)^{-1}\widehat{\bfW}_r^T\|\,\|\cbfB(s)\|}
{\|\hatcbfH_r(s)\|}\\
\mbox{cond}_{\subcbfC}(\cbfH_r(s))&=\frac{\|\cbfC(s)\|\,\|\bfV_r\cbfK_r(s)^{-1}\cbfB_r(s)\|}{\|\cbfH_r(s)\|}
\end{align*} 
   \end{theorem}
   \textsc{Proof:} Note that for all $s\in\IC$ for which $\cbfH_r$ and $\hatcbfH_r$ are both analytic, 
   {\small
   \begin{align*} 
    \|\cbfH_r(s)-\hatcbfH_r(s)\|&=\|\cbfC(s)\left(\bfV_r\cbfK_{r}(s)^{-1}\bfW_r^T - \widehat{\bfV}_r\hatcbfK_{r}(s)^{-1}\widehat{\bfW}_r^T\right)\cbfB(s)\|\\
    & =\|\cbfC(s)\left(\cbfQ_{r}(s)-\hatcbfQ_{r}(s)\right)\cbfK(s)^{-1}\cbfB(s)\|\\
    & =\|\cbfC(s)\left(\left(\bfI-\hatcbfQ_{r}(s)\right)\cbfQ_{r}(s)-\hatcbfQ_{r}(s)\left(\bfI-\cbfQ_{r}(s)\right)\right)\cbfK(s)^{-1}\cbfB(s)\|
\end{align*} 
So, 
\begin{align*} 
   \|\cbfH_r(s)-\hatcbfH_r(s)\|  & \leq\|\cbfC(s)\left(\bfI-\hatcbfQ_{r}(s)\right)\cbfQ_{r}(s)\cbfK(s)^{-1}\cbfB(s)\| \\
    &\qquad\qquad +\|\cbfC(s)\hatcbfQ_{r}(s)\left(\bfI-\cbfQ_{r}(s)\right)\cbfK(s)^{-1}\cbfB(s)\|\\
   & \leq\|\cbfC(s)\left(\bfI-\hatcbfQ_{r}(s)\right)\cbfQ_{r}(s)\cbfK(s)^{-1}\cbfB(s)\|\\
    &\qquad\qquad +\|\cbfC(s)\cbfK(s)^{-1}\hatcbfP_{r}(s)\left(\bfI-\cbfP_{r}(s)\right)\cbfB(s)\|\\
    & \leq\|\cbfC(s)\left(\bfI-\hatcbfQ_{r}(s)\right)\,\left(\bfI -  \widehat{\bfXi} \right) \bfXi\,\cbfQ_{r}(s)\cbfK(s)^{-1}\cbfB(s)\|\\
    &\qquad\qquad +\|\cbfC(s)\cbfK(s)^{-1}\hatcbfP_{r}(s)\,\widehat{\bfPi}\left(\bfI - \bfPi \right)\,\left(\bfI-\cbfP_{r}(s)\right)\cbfB(s)\|\\
    & \leq\|\cbfC(s)\| \left\|\bfI-\hatcbfQ_{r}(s)\right\|\,\left\|\left(\bfI -  \widehat{\bfXi} \right) \bfXi\right\|\,\|\cbfQ_{r}(s)\cbfK(s)^{-1}\cbfB(s)\|\\
    &\qquad\qquad +\|\cbfC(s)\cbfK(s)^{-1}\hatcbfP_{r}(s)\|\,\|\widehat{\bfPi}\left(\bfI - \bfPi \right)\|\,\|\bfI-\cbfP_{r}(s)\|\,\|\cbfB(s)\|\\
    & \leq\|\cbfC(s)\| \frac{\sin\Theta(\widehat{\mathcal{V}}_{r},\,\mathcal{V}_r) }{\cos\Theta(\widehat{\mathfrak{Q}}_{r}(s),\,\widehat{\mathcal{V}}_r)} 
    \|\cbfQ_{r}(s)\cbfK(s)^{-1}\cbfB(s)\|\\
    &\qquad\qquad +\|\cbfC(s)\cbfK(s)^{-1}\hatcbfP_{r}(s)\|\,
    \frac{\sin\Theta(\widehat{\mathcal{W}}_{r},\,\mathcal{W}_r) }{\cos\Theta(\mathfrak{P}_{r}(s),\,\mathcal{W}_r)}\,\|\cbfB(s)\|
   \end{align*} 
   }
   \begin{align*} 
    \qquad \qquad & \leq\mbox{cond}_{\subcbfC}(\cbfH_r(s))
     \frac{\sin\Theta(\widehat{\mathcal{V}}_{r},\,\mathcal{V}_r) }{\cos\Theta(\widehat{\mathfrak{Q}}_{r}(s),\,\widehat{\mathcal{V}}_r)} 
    \|\cbfH_{r}(s)\|\\
    &\qquad\qquad +\|\hatcbfH_{r}(s)\|\,
    \frac{\sin\Theta(\widehat{\mathcal{W}}_{r},\,\mathcal{W}_r) }{\cos\Theta(\mathfrak{P}_{r}(s),\,\mathcal{W}_r)}\,
   \mbox{cond}_{\subcbfB}(\hatcbfH_r(s))
   \end{align*} 
Maximizing over $s=\imath\omega$ with $\omega\in\IR$ gives
  {\small
 \begin{align*} 
  \|\cbfH_r-\hatcbfH_r\|_{\hinf}
  & \leq \max_{\omega\in\IR} \mbox{cond}_{\subcbfC}(\cbfH_r(\imath\omega))
     \frac{\sin\Theta(\widehat{\mathcal{V}}_{r},\,\mathcal{V}_r) }{\min_{\omega\in\IR}\cos\Theta(\widehat{\mathfrak{Q}}_{r}(s),\,\widehat{\mathcal{V}}_r)} 
    \|\cbfH_{r}\|_{\hinf}\\
    &\qquad\quad +
         \max_{\omega\in\IR} \mbox{cond}_{\subcbfB}(\hatcbfH_r(\imath\omega))
    \frac{\sin\Theta(\widehat{\mathcal{W}}_{r},\,\mathcal{W}_r) }{\min_{\omega\in\IR}\cos\Theta(\mathfrak{P}_{r}(s),\,\mathcal{W}_r)}\,
    \|\hatcbfH_{r}\|_{\hinf}
  \end{align*} 
  } 
which leads immediately to the conclusion.  $\Box$

\subsection{Illustrative examples} \label{examples}
The process to be modeled arises in cooling within
 a rolling mill and is modeled as
boundary control of a two dimensional heat equation. A finite
element discretization
 results in a descriptor system of the form
$$ \bfE \dot \bfx(t) =  \bfA \bfx(t) + \bfB \bfu(t), ~~~~
              y(t)  =  \bfC \bfx(t).$$
 where $\bfA,\bfE \in \IR^{5177
    \times 5177}$, $\bfB \in \IR^{5177 \times 7} $, $\mathbf{C}
\in \IR^{6 \times
    5177}$. For simplicity, we focus
on a SISO full-order subsystem that relates the sixth
          input  to the second output.
      For details regarding the modeling, discretization,
      optimal control design,  and model reduction,  see \cite{benner2004slc,benner2004ens}.

We show the results of interpolatory model reduction using an \emph{ad hoc} choice of interpolation points:  $6$ logarithmically spaced points between $10^{0.5}$ and $10$; and an $\h2$-optimal choice of interpolation points obtained by the  method of \cite{gugercin2008hmr}. For each case, we reduce the system order to $r=6$ using first \emph{exact} interpolatory model reduction (i.e., the linear systems are solved directly) and then with inexact model reduction with varying choices of termination criteria. 
The resulting reduced-order models are denoted by $\cbfH_r(s)$ and $\hatcbfH_r(s)$, respectively. 
To see the effect of 
the choice of interpolation points on the underlying model reduction problem,
we vary the relative residual termination tolerance, $\varepsilon$
 between $10^{-1}$ and $10^{-10}$ and show how quickly 
$\hatcbfH_r(s)$ converges to $\cbfH_r(s)$ for both the \emph{ad hoc} selection and the $\h2$-optimal selection of interpolation points. 
 Table \ref{table:HminusHr} shows the relative 
$\hinf$ error between $\hatcbfH_r(s)$ and $\cbfH_r(s)$ as $\varepsilon$ decreases. 
 For the $\h2$-optimal choice of interpolation points,
$\hatcbfH_r(s)$ converges to $\cbfH_r(s)$ as $\varepsilon$ decreases, for the \emph{ad hoc} choice of points, there is almost no improvement
in accuracy until $\varepsilon = 1\times 10^{-6}$.
 \begin{table} [ht]
\begin{center}
 \begin{tabular}{c|*{3}{c}}
      $\varepsilon$ & $\h2$-optimal $\{\sigma_i \}$ & \emph{ad hoc} $\{\sigma_i \}$    \\       
      \hline
     $ 10^{-1}$ &$7.22 \times 10^{-1}$& $5.05 \times 10^{-1}$  \\
      $ 10^{-2}$   & $2.00\times 10^{-1}$  & $1.64 \times 10^{-1}$  \\
  $ 10^{-3}$   & $4.27 \times 10^{-2}$ & $4.11 \times 10^{-1}$ \\
  $ 10^{-4}  $ & $1.07 \times 10^{-2}$ & $ 2.38 \times 10^{-1} $ \\
  $ 10^{-5}  $ & $2.76 \times 10^{-4} $ & $5.62 \times 10^{-1} $ \\
  $ 10^{-6}  $ & $2.56 \times 10^{-5}$ & $2.13\times 10^{-2}$ \\
  $ 10^{-7} $  & $2.91\times 10^{-6} $& $3.52 \times 10^{-3}$  \\
  $ 10^{-8} $  & $1.51 \times 10^{-7}$ & $6.18 \times 10^{-5} $ \\
 $ 10^{-9}$ & $2.07 \times 10^{-8}$ & $1.76 \times 10^{-5}$  \\
  $ 10^{-10}$   & $2.17 \times 10^{-9}$ & $5.15 \times 10^{-6} $ \\
  \end{tabular}
   \end{center}
       \caption{The relative error $\displaystyle \frac{\bigl\|\cbfH_{r}-\hatcbfH_{r}\bigr\|_{\mathcal{H}_{\infty}}}{
        \bigl\|\cbfH_{r}\bigr \|_{\mathcal{H}_{\infty}}}$ as $\varepsilon$ varies}
   \label{table:HminusHr}
\end{table}

The behavior exhibited in Table \ref{table:HminusHr} becomes clearer once we inspect the subspace angles between the exact
interpolatory subspaces $\mathcal{V}_r$, $\mathcal{W}_r$  and the inexact ones $\widehat{\mathcal{V}}_{r}$
and $\widehat{\mathcal{W}}_{r}$. Table \ref{table:angles} shows the sine of the angle 
between the exact and inexact interpolatory subspaces as $\varepsilon$ varies. 
While the gap 
decreases significantly as $\varepsilon$ decreases  for an $\h2$-optimal selection of interpolation points, there is a much smaller improvement 
in the gap with respect to $\varepsilon$ for an \emph{ad hoc} choice of points. This behavior will be re-visited in 
more detail in \S \ref{sec:condnumber} revealing that the $\h2$-optimal (or good) interpolation points are expected to produce reduced order models that are more robust with respect to perturbations due to inexact solves. 
{
\begin{table} [ht]
\begin{center}
 \begin{tabular}{c||c|c||c|c|}
& \multicolumn{2}{c||}{{ $\sin\Theta({\mathcal{V}}_r,\widehat{\mathcal{V}}_{r})$ }} & \multicolumn{2}{c|}{{ $\sin\Theta({\mathcal{W}}_r,\widehat{\mathcal{W}}_{r})$ }} \\
  $\varepsilon$ &  $\h2$-optimal $\{\sigma_i\}$  & \emph{ad hoc} $\{\sigma_i\}$ &  $\h2$-optimal $\{\sigma_i\}$  & 
  \emph{ad hoc} $\{\sigma_i\}$    \\
\hline
 $10^{-1}$ & $\relnorm{9.85}{1}$ & $\relnorm{9.99}{1}$ 
 & $\relnorm{9.99}{1}$ & $\relnorm{9.99}{1}$ \\
 $10^{-2}$   &  $\relnorm{1.99}{1}$ & $\relnorm{9.99}{1}$ 
& $\relnorm{9.97}{1}$ & $\relnorm{9.93}{1}$  \\
 $10^{-3}$    &  $\relnorm{2.36}{2}$ & $\relnorm{9.99}{1}$
 & $\relnorm{4.87}{1}$ & $\relnorm{9.83}{1}$ \\
 $10^{-4}$ 	  &  $\relnorm{4.39}{3}$ & $\relnorm{9.60}{1}$  
 & $\relnorm{6.38}{2}$ & $\relnorm{9.99}{1}$
 \\
 $10^{-5}$ & $\relnorm{2.72}{4}$ & $ \relnorm{5.80}{1}$  
  & $\relnorm{7.09}{3}$ & $\relnorm{7.20}{1}$ \\
$10^{-6}$ & $\relnorm{2.90}{5}$ & $\relnorm{4.57}{2}$  
 & $\relnorm{9.88}{4}$ & $\relnorm{1.19}{1}$\\
$10^{-7}$ & $\relnorm{3.46}{6}$ &  $\relnorm{6.90}{3}$ 
 & $\relnorm{6.87}{5}$ & $\relnorm{2.00}{2}$  \\
$10^{-8}$ & $\relnorm{3.85}{7}$ &  $\relnorm{7.92}{4}$
 & $\relnorm{6.71}{6}$ & $\relnorm{2.26}{3}$\\
$10^{-9}$ & $\relnorm{3.63}{8}$ & $\relnorm{1.01}{4}$  
 & $\relnorm{9.16}{7}$ & $\relnorm{2.60}{4}$\\
$10^{-10}$ &  $\relnorm{2.71}{9}$ & $\relnorm{1.28}{5}$
 & $\relnorm{6.35}{8}$ & $\relnorm{3.10}{5}$ \\
 \end{tabular}
 \end{center}
 \caption{$r = 6$; $\sin\Theta({\mathcal{V}}_r,\widehat{\mathcal{V}}_{r})$
 and $\sin\Theta({\mathcal{W}}_r,\widehat{\mathcal{V}}_{r})$ as $\varepsilon$ varies}
 \label{table:angles}
 \end{table}
 }
 
\section{Backward error} \label{sec:backwarderror}
Instead of seeking bounds on how much an inexactly computed reduced model differs from an exactly computed counterpart, one may view an inexactly computed reduced order model as an exactly computed reduced order model of a perturbed full order system.   That is, we wish to find a full order system
\begin{equation} \label{backerrSys}
\tildecbfH(s)= \tildecbfC(s)\tildecbfK(s)^{-1}\tildecbfB(s)
\end{equation}
so that the inexactly computed reduced model for $\cbfH(s)= \cbfC(s)\cbfK(s)^{-1}\cbfB(s)$ 
would be an \emph{exactly} computed interpolatory reduced model for $\tildecbfH(s)$.   Given left and right tangential interpolation data as in 
(\ref{lefttang}) and (\ref{righttang}) that has contributed toward producing the inexactly computed interpolatory reduced model $\hatcbfH_r(s)$, 
find $\tildecbfH(s)$ as in (\ref{backerrSys}) so that 
\begin{align*}
\bfc_i^T \tildecbfH(\mu_i)&=\bfc_i^T \hatcbfH_r(\mu_i)\quad\mbox{ for }i=1,\,\ldots,\,r, \mbox{ and }\\
 \tildecbfH(\sigma_j) & \bfb_j = \hatcbfH_r(\sigma_j) \bfb_j \quad \mbox{ for }j=1,\,\ldots,\,r.
\end{align*}
and so that $\hatcbfH_r$ could have been computed from the perturbed system $\tildecbfH$ from the given tangential interpolation data via an exact computation.
Specifically, given computed (inexact) projecting bases 
\begin{equation*}   
\widehat{\bfV}_r =\  \left[~ \widehat{\bfv}_1,~\cdots,~ \widehat{\bfv}_r~\right]
\qquad
\widehat{\bfW}_r^T =
\ \left[\begin{array}{c}
\widehat{\bfw}_1^T\\
\vdots \\
\widehat{\bfw}_r^T
\end{array}\right].
\end{equation*}
as in (\ref{defineVhatWhat}), and a resulting (inexact) reduced order coprime realization 
$$
\hatcbfH_r(s)= \hatcbfC_r(s) \hatcbfK_r(s)^{-1} \hatcbfB_r(s),
$$ 
find a full-order system $\tildecbfH(s)= \tildecbfC(s)\tildecbfK(s)^{-1}\tildecbfB(s)$ so that left and right interpolation conditions hold:
\begin{align}
\bfc_i^T\tildecbfC(\mu_i)&=\widehat{\bfw}_i^T\tildecbfK(\mu_i) \quad\mbox{ for }i=1,\,\ldots,\,r, \label{backerrleftInt}\\
\tildecbfK(\sigma_j)\widehat{\bfv}_j&=\tildecbfB(\sigma_j) \bfb_j \quad \mbox{ for }j=1,\,\ldots,\,r, \label{backerrrightInt}
\end{align}
and so that
\begin{equation}  \label{backerrSysRed}
\hatcbfK_r(s)= \widehat{\bfW}_r^{T} \tildecbfK(s) \widehat{\bfV}_r,\quad
\hatcbfB_r(s) = \widehat{\bfW}_r^{T}\tildecbfB(s), \quad 
\mbox{and}\quad \hatcbfC_r(s) = \tildecbfC(s) \widehat{\bfV}_r,
\end{equation}
There (typically) will be an infinite number of possible systems, $\tildecbfH$, that are consistent with the computed reduced system $\hatcbfH_r$ in this sense --- we are interested in those that are \emph{close} to the original system $\cbfH$ with respect to a convenient system norm such as $\hinf$ or $\h2$.  In order to proceed, it is convenient to restrict the class of backwardly compatible systems, $\tildecbfH$.  We consider those that have realizations that are \emph{constant} perturbations from the corresponding original system factors:
\begin{equation}  \label{backPertForm}
\tildecbfK(s)= \cbfK(s) +\bfF, \quad
\tildecbfB(s) =  \cbfB(s) +\bfE, \quad 
\mbox{and}\quad \tildecbfC(s) =  \cbfC(s) +\bfG.
\end{equation}
where $\bfE$, $\bfF$, and $\bfG$ are \emph{constant} matrices. 
The conditions (\ref{backerrleftInt}), (\ref{backerrrightInt}), and (\ref{backerrSysRed}) impose constraints on $\bfE$, $\bfF$, and $\bfG$.  
Indeed, (\ref{backerrleftInt}) and (\ref{backerrrightInt}) imply that 
\begin{align*}
\widehat{\bfw}_i^T\bfF+\bfxi_i^T&=\bfc_i^T\bfE \quad\mbox{ for }i=1,\,\ldots,\,r, \mbox{ and}\\
\bfF\widehat{\bfv}_j+\bfeta_j&=\bfG \bfb_j \quad \mbox{ for }j=1,\,\ldots,\,r.
\end{align*}
(\ref{backerrSysRed}) implies that 
$$
\widehat{\bfW}_r^{T} \bfF\widehat{\bfV}_r=\mathbf{0},\quad
\widehat{\bfW}_r^{T} \bfG=\mathbf{0},\quad\mbox{and}\quad
 \bfE\widehat{\bfV}_r=\mathbf{0}.
$$
Taken together, we find that backward perturbations of the form (\ref{backPertForm}) can exist only if 
\begin{equation} \label{backerrCompat}
\bfxi_i^T\widehat{\bfV}_r=\mathbf{0} \quad\mbox{ for }i=1,\,\ldots,\,r, \mbox{ and }
\widehat{\bfW}_r^{T}\bfeta_j= \mathbf{0}\quad \mbox{ for }j=1,\,\ldots,\,r.
\end{equation}
Thus, we find constraints on the inexact interpolation residuals $\bfxi_i$ and $\bfeta_j$ in order for a backwardly compatible 
system of the form (\ref{backPertForm}) to exist.  More complicated perturbation classes than (\ref{backPertForm}) may be
considered that would allow us to remove the conditions (\ref{backerrCompat}), of course, but instead we 
choose to focus on a computational framework that guarantees (\ref{backerrCompat}). 
The Biconjugate Gradient Algorithm (BiCG) will be an example of an iterative solution 
strategy that fits this framework \cite{Kap09,barrettTemplates1994}; others can be constructed without difficulty, although many standard strategies such as GMRES, do not fit this framework.  


\subsection{The Petrov-Galerkin Framework for Inexact Solves} \label{sec:pg}
We have observed above that (\ref{backerrCompat}) is necessary for there to be a well-defined backward error of the form (\ref{backPertForm})  to exist. 
The simplest framework within which one may generate reduced order models that are guaranteed to satisfy this condition involves a Petrov-Galerkin 
formalism for producing approximate solutions to (\ref{lefttangkrylov}) and (\ref{righttangkrylov}).  For simplicity, we restrict our discussion to the case that $\mu_i=\sigma_i$ (identical left and right interpolation points).

Let 
$\mathcal{P}_{N}$ and $\mathcal{Q}_{N}$ be $N$-dimensional subspaces of 
 $\mathbb{C}^{n}$ satisfying a nondegeneracy condition:
 $\left(\cbfK(\sigma_i) \mathcal{P}_{N}\right)^{\perp}\cap \mathcal{Q}_{N}= \{0\}$ 
 for all shifts, $\sigma_i$ to be considered. 
 The \emph{Petrov-Galerkin framework} for generating approximate solutions to the interpolation conditions (\ref{lefttangkrylov}) and (\ref{righttangkrylov}) proceeds as follows:
 \begin{align} \label{eqn:pgframework}
 \mbox{Find }&\widetilde{\mathbf{v}}_{j} \in \mathcal{P}_{N} \mbox{ so that }
 \cbfK(\sigma_{j}) \widetilde{\bfv}_j - \cbfB(\sigma_j)\bfb_j\perp  \mathcal{Q}_{N} \quad\mbox{and} \nonumber \\
  &\mbox{find }\widetilde{\mathbf{w}}_{j} \in \mathcal{Q}_{N} \mbox{ so that }
   \cbfK(\sigma_{j})^T\widetilde{\bfw}_j -  \cbfC(\sigma_j)^T\bfc_i  \perp \mathcal{P}_{N}
 \end{align} 
 Computed quantities generated within a Petrov-Galerkin framework will be denoted with a ``tilde'' to distinguish them from earlier ``hat" quantities  where no structure was assumed in the inexact solves. 
 The following theorem asserts that if a reduced order model is computed within a Petrov-Galerkin framework (\ref{eqn:pgframework}),
 then one can obtain a structured backward error that throws the effect of inexact solves back onto a perturbation on the original dynamical system. 
\begin{theorem} \label{backwardHrtilde}
Given a full order model $\cbfH(s)= \cbfC(s)\cbfK(s)^{-1}\cbfB(s)$, interpolation points $\{\sigma_j\}_{j=1}^r$, and tangent directions
$\{ \bfb_i\}_{i=1}^r$ and $\{ \bfc_i\}_{i=1}^r$, 
 let the inexact solutions 
$\widetilde{\bfv}_j$ for $\cbfK(\sigma_j)^{-1}\cbfB(\sigma_j)\bfb_j$ and $\widetilde{\bfw}_j$ for $\cbfK(\sigma_j)^{-T}\cbfC(\sigma_j)^T\bfc_j$ 
be obtained in a Petrov-Galerkin framework as in (\ref{eqn:pgframework}).
Let $\widetilde{\bfV}_r$ and $\widetilde{\bfW}_r$
denote the corresponding inexact interpolatory bases; i.e.
\begin{eqnarray}  \label{VrWrtildebases}
\widetilde{\bfV}_r = \left[~ \widetilde{\bfv}_1,~~\cdots,~~ \widetilde{\bfv}_r~\right] ~~~~
{\rm and}~~~~ \widetilde{\bfW}_r = \left[~ \widetilde{\bfw}_1,~~\cdots,~~ \widetilde{\bfw}_r~\right].
\end{eqnarray}

Define residuals 
$$
{\bfeta}_j  = \cbfK(\sigma_j) \widetilde{\bfv}_j-\cbfB(\sigma_j)\bfb_j ~~~~~~{\rm and}~~~~~~{\bfxi}_j = 
\cbfK(\sigma_j)^T \widetilde{\bfw}_j-\cbfC(\sigma_j)^T\bfc_j,
$$
residual matrices 
\begin{equation} \label{eqn:RbRc}
\mathbf{R}_{\mathbf{b}}=\left[{\bfeta}_{1},\,{\bfeta}_{2},\,
\ldots,{\bfeta}_{r}\right], \qquad
\mathbf{R}_{\mathbf{c}}=\left[{\bfxi}_{1},\,{\bfxi}_{2},\,
\ldots,{\bfxi}_{r}\right],
\end{equation}
and the rank $2r$ matrix
\begin{equation} \label{E2r}
\mathbf{F}_{2r}=\mathbf{R}_{\mathbf{b}}
(\widetilde{\mathbf{W}}_{r}^{T} \widetilde{\mathbf{V}}_{r})^{-1}
\widetilde{\mathbf{W}}_{r}^{T}+ \widetilde{\mathbf{V}}_{r}
(\widetilde{\mathbf{W}}_{r}^{T} \widetilde{\mathbf{V}}_{r})^{-1}
\mathbf{R}_{\mathbf{c}}^{T}.
\end{equation}
Let $\tildecbfH_r(s) = \tildecbfC_r(s)\tildecbfK_r(s)^{-1}\tildecbfB_r(s)$ denote the computed inexact reduced model via the Petrov-Galerkin process where
\begin{equation} \label{Hrtilde}
 \tildecbfK_r(s) = \widetilde{\mathbf{W}}_{r}^{T}
\cbfK(s) \widetilde{\mathbf{V}}_{r},~~
\tildecbfB_r(s) = \widetilde{\mathbf{W}}_{r}^{T}
\cbfB(s), \mbox{ ~~~and~~~ }
\tildecbfC_r(s) =
\cbfC(s) \widetilde{\mathbf{V}}_{r}.
\end{equation}
Then, $\tildecbfH_r(s)$ exactly tangentially interpolates the perturbed full-order model
\begin{equation}
\tildecbfH(s) = \cbfC(s)(\cbfK(s)+\bfF_{2r})^{-1}
\cbfB(s),
\end{equation}
at each $\sigma_i$: 
\begin{align*}
\tildecbfH(\sigma_i)\bfb_i &= \tildecbfH_r(\sigma_i)\bfb_i, \quad
\bfc_i^T \tildecbfH(\sigma_i)= \bfc_i^T\tildecbfH_r(\sigma_i), \\
\mbox{and}&\quad
\bfc_i^T \tildecbfH'(\sigma_i) \bfb_i= \bfc_i^T \tildecbfH_r'(\sigma_i)\bfb_i\quad \mbox{for each }i=1,\ldots,r.
\end{align*}
\end{theorem}
\textsc{Proof:}  The computed model, $\tildecbfH_r(s)$, will (exactly) tangentially interpolate a perturbed model 
$\tildecbfH(s) = \cbfC(s)(\cbfK(s)+\bfF)^{-1}\cbfB(s)$ provided the following interpolation conditions hold:
$$
\left(\cbfK(\sigma_i)+\bfF\right) \widetilde{\bfv}_i =\cbfB(\sigma_i)\bfb_i \ \mbox{and} \ \widetilde{\bfw}_i^T\left(\cbfK(\sigma_i)+\bfF\right)= \bfc_i^T\cbfC(\sigma_i)  \quad \mbox{for }i=1,\ldots,r.
$$
Equivalently, these can be interpreted as conditions on the perturbation $\bfF$.  Rewriting this using notation defined above, $\bfF$ must satisfy
\begin{equation}\label{backpertcond}
\bfF \widetilde{\mathbf{V}}_{r} = \mathbf{R}_{\mathbf{b}} \quad \mbox{and}\quad  \widetilde{\mathbf{W}}_{r}^{T} \bfF=\mathbf{R}_{\mathbf{c}}^{T}.
\end{equation}

The Petrov-Galerkin framework guarantees  $\widetilde{\mathbf{W}}_{r}^{T}\mathbf{R}_{\mathbf{b}} = \mathbf{0}$ and 
$\mathbf{R}_{\mathbf{c}}^{T}\widetilde{\mathbf{V}}_{r}= \mathbf{0}$.  Substitution of $\bfF_{2r}$ from (\ref{E2r}) into (\ref{backpertcond}) verifies that
$\bfF_{2r}$ is a perturbation to $\cbfK(s)$ for which the computed (inexact) vectors become (exact) interpolation vectors.

Note that since $\widetilde{\mathbf{W}}_{r}^{T}\mathbf{F}_{2r}\widetilde{\mathbf{V}}_{r}=\mathbf{0}$,
$$
\tildecbfK_r(s) =  \widetilde{\mathbf{W}}_{r}^{T}
\cbfK(s) \widetilde{\mathbf{V}}_{r} = \widetilde{\mathbf{W}}_{r}^{T}
( \cbfK(s) + \bfF_{2r})\widetilde{\mathbf{V}}_{r}.
$$
Consequently, 
the reduced model
$\tildecbfH_r(s)$ obtained by inexact solves in (\ref{Hrtilde}) is what one would have 
obtained by exact interpolatory model reduction of $\tildecbfH(s)$. $\Box$\\

\begin{theorem} \label{thm:F2r_fro_norm}
Assume the hypotheses of Theorem \ref{backwardHrtilde} 
and that $\widetilde{\bfW}_r^T \widetilde{\bfV}_r$ is nonsingular.  Define an oblique projector,   $ {\widetilde{\bfPhi}}_r= \widetilde{\bfV}_r(\widetilde{\bfW}_r^T \widetilde{\bfV}_r)^{-1} \widetilde{\bfW}_r^T$.  
The backward perturbation $\bfF_{2r}$ given in Theorem \ref{backwardHrtilde} satisfies 
$$
\|\mathbf{F}_{2r}\|_{F}\leq \sqrt{r}\,\|{\widetilde{\bfPhi}}_r\|\cdot \left(
\max_{i} \frac{\|\bfeta_i\|} {\|\widetilde{\bfv}_i\|}
\varsigma_{\min}(\widetilde{\mathbf{V}}_{r}\bfD)^{-1}+
\max_{i} \frac{\|\bfxi_i\|} {\|\widetilde{\bfw}_i\|} 
\varsigma_{\min}(\widetilde{\mathbf{W}}_{r}\bfD)^{-1}\right)
$$
where $\varsigma_{\min}$ denotes the smallest singular value and 
$\|\bfM\|_F=\sqrt{\mbox{\textsf{trace}}(\bfM^T\bfM)}$ denotes the Frobenius norm of a matrix, $\bfM$.
\end{theorem}
\textsc{Proof:} 
Note that 
$$
\|\mathbf{F}_{2r}\|_{F}\leq
\|\mathbf{R}_{\mathbf{b}}(\widetilde{\mathbf{W}}_{r}^{T} \widetilde{\mathbf{V}}_{r})^{-1}
\widetilde{\mathbf{W}}_{r}^{T}\|_F +
\|\widetilde{\mathbf{V}}_{r}(\widetilde{\mathbf{W}}_{r}^{T} \widetilde{\mathbf{V}}_{r})^{-1}\mathbf{R}_{\mathbf{c}}^{T}\|_F.
$$
Let $\widetilde{\mathbf{V}}_{r}$ have an orthogonal factorization as $\widetilde{\mathbf{V}}_{r}=\bfQ_v\bfL_v$ with $\bfQ_v^*\bfQ_v=\bfI$.  Then
\begin{align*}
\|\mathbf{R}_{\mathbf{b}}(\widetilde{\mathbf{W}}_{r}^{T} \widetilde{\mathbf{V}}_{r})^{-1}\widetilde{\mathbf{W}}_{r}^{T}\|_F
& = \|\mathbf{R}_{\mathbf{b}} \bfL_v^{-1} \bfL_v(\widetilde{\mathbf{W}}_{r}^{T} \widetilde{\mathbf{V}}_{r})^{-1}\widetilde{\mathbf{W}}_{r}^{T}\|_F\\
&\leq \|\mathbf{R}_{\mathbf{b}} \bfL_v^{-1}\|_F \cdot \|\bfL_v(\widetilde{\mathbf{W}}_{r}^{T} \widetilde{\mathbf{V}}_{r})^{-1}\widetilde{\mathbf{W}}_{r}^{T}\|\\
&\leq \|\mathbf{R}_{\mathbf{b}} \bfL_v^{-1}\|_F \cdot \|{\widetilde{\bfPhi}}_r\|\\
&\leq \|\mathbf{R}_{\mathbf{b}}\widetilde{\bfD}_{v}(\bfL_v\widetilde{\bfD}_{v})^{-1}\|_F \cdot \|{\widetilde{\bfPhi}}_r\|\\
&\leq \|\mathbf{R}_{\mathbf{b}}\widetilde{\bfD}_{v}\|_F\cdot\|(\bfL_v\widetilde{\bfD}_{v})^{-1}\| \cdot \|{\widetilde{\bfPhi}}_r\|
\end{align*}
where we have introduced a diagonal scaling matrix 
$$\widetilde{\bfD}_{v}=\mbox{diag}(1/\|\widetilde{\bfv}_1\|,\,1/\|\widetilde{\bfv}_2\|,\,\ldots,\,1/\|\widetilde{\bfv}_r\|).$$
Easily one sees $\displaystyle \|\mathbf{R}_{\mathbf{b}}\widetilde{\bfD}_{v}\|_F\leq \sqrt{r}\max_i\frac{\|\bfeta_i\|}{\|\widetilde{\bfv}_i\|}$.  For the remaining term, 
note that 
$$
\|(\bfL_v\widetilde{\bfD}_{v})^{-1}\|=\left(\min_{\bfx}\frac{\|\widetilde{\mathbf{V}}_{r}\widetilde{\bfD}_{v}\bfx\|}{\|\bfx\|}\right)^{-1}
= \varsigma_{\min}\!\left(\widetilde{\mathbf{V}}_{r}\widetilde{\bfD}_{v}\right)^{-1}
$$
A similar bound for 
{\small $\|\widetilde{\mathbf{V}}_{r}(\widetilde{\mathbf{W}}_{r}^{T} \widetilde{\mathbf{V}}_{r})^{-1}\mathbf{R}_{\mathbf{c}}^{T}\|_F$} is produced by an analogous process, which leads then to the final estimate for $\|\mathbf{F}_{2r}\|_{F}$.
 $\Box$

Note that 
the perturbation $\bfF_{2r}$ is completely determined by accessible, computed quantities. Hence, one can use
$\bfF_{2r}$  to determine how accurately one must solve the underlying linear systems in order to assure system fidelity of a given order.

\begin{theorem} \label{thm:Fper_bound}
If  {\small $\|\bfF_{2r}\| <1/\| \cbfK(s)^{-1}\|_{\hinf}$} then
{\small
$$
\|\cbfH(s) -  \tildecbfH(s)\|_{\h2}\leq
\frac{ \|\cbfC(s)\cbfK(s)^{-1}\|_{\h2}\ \|\cbfK(s)^{-1}\cbfB(s)\|_{\hinf}}
{1-\| \cbfK(s)^{-1}\|_{\hinf}\,\|\bfF_{2r}\| } \|\bfF_{2r}\|
$$
}
\end{theorem}

\textsc{Proof:} 
The system-wise backward error associated with inexact solves may be written as
\begin{align*}
\cbfH(s) -  \tildecbfH(s) &= \cbfC(s)\cbfK(s)^{-1}\cbfB(s) -  \cbfC(s)\left(\cbfK(s)+\bfF_{2r}\right)^{-1}\cbfB(s) \nonumber \\ 
& =  \cbfC(s)\cbfK(s)^{-1}\bfF_{2r}\left(\cbfK(s)+\bfF_{2r}\right)^{-1}\cbfB(s) \\
& =  \cbfC(s)\cbfK(s)^{-1}\bfF_{2r}\left(\bfI+\cbfK(s)^{-1}\bfF_{2r}\right)^{-1}\cbfK(s)^{-1}\cbfB(s) \nonumber
\end{align*}
Define $\cbfM(s)=\bfF_{2r}\left(\bfI+\cbfK(s)^{-1}\bfF_{2r}\right)^{-1}$ and observe that
\begin{align*}
&\|\cbfH(s) -  \tildecbfH(s)\|_{\h2}^2=\frac{1}{2\pi}\int_{-\infty}^{\infty} \|\cbfC(\imath \omega)\cbfK(\imath \omega)^{-1}\,\cbfM(\imath \omega)\,
\cbfK(\imath \omega)^{-1}\cbfB(\imath \omega)\|_F^2\ d\omega\\
&\qquad \leq \frac{1}{2\pi}\int_{-\infty}^{\infty} \|\cbfC(\imath \omega)\cbfK(\imath \omega)^{-1}\|_F^2\,\cdot\,\|\,\cbfM(\imath \omega)\, \|^2\, 
\cdot\, \|\cbfK(\imath \omega)^{-1}\cbfB(\imath \omega)\|^2\ d\omega\\
& \quad\leq \left(\frac{1}{2\pi}\int_{-\infty}^{\infty} \|\cbfC(\imath \omega)\cbfK(\imath \omega)^{-1}\|_F^2\,d\omega\right)\,\cdot\, \max_\omega\| \cbfM(\imath \omega)\ \|^2\, \cdot\, 
\max_\omega\| \cbfK(\imath \omega)^{-1}\cbfB(\imath \omega)\|^2\\
&\quad \leq \,\|\cbfC(s)\cbfK(s)^{-1}\|_{\h2}^2\,\cdot\, \|\cbfK(s)^{-1}\cbfB(s)\|_{\hinf}^2 \, \cdot\, 
\| \cbfM(s)\|_{\hinf}^2.\end{align*}
To estimate $\| \cbfM(s)\|_{\hinf}$, a rearrangement of the definition of $\cbfM(s)$ provides
$$
\cbfM(s)= \left(\bfI-\cbfM(s)\cbfK(s)^{-1}\right)\bfF_{2r}.
$$
So we have immediately,
\begin{align*}
\| \cbfM(s)\|_{\hinf}=& \max_{\omega\in\IR} \|\cbfM(\imath \omega)\|\leq 
\max_{\omega\in\IR} \|\bfI-\cbfM(\imath \omega)\cbfK(\imath \omega)^{-1}\|\, \cdot\, \|\bfF_{2r}\|\\
& \leq  \left(1 + \max_{\omega\in\IR} \|\cbfM(\imath \omega)\cbfK(\imath \omega)^{-1}\|\right)\, \|\bfF_{2r}\|\\
& \leq  \left(1 + \| \cbfM(s)\|_{\hinf}\, \|\cbfK(s)^{-1}\|_{\hinf}\right)\, \|\bfF_{2r}\|
\end{align*}
Since $\| \cbfK(s)^{-1}\|_{\hinf}\,\|\bfF_{2r}\| <1$, this last expression can be rearranged to obtain
$$
\| \cbfM(s)\|_{\hinf} \leq\frac{\|\bfF_{2r}\|}
{1-\| \cbfK(s)^{-1}\|_{\hinf}\,\|\bfF_{2r}\| }
$$
which implies the conclusion. $\Box $

%
%
By combining Theorem \ref{globalFwdBound} with Theorem \ref{thm:subspacebound} or  combining Theorem \ref{thm:F2r_fro_norm} with Theorem \ref{thm:Fper_bound}, we  approach our goal of connecting quantities that we have control over, such as the termination threshold, $\varepsilon$,  
to relevant system theoretic errors, $\|\cbfH_r-\hatcbfH_r\|$ and $\| \cbfH -  \tildecbfH\|$, which are quantities we would like to control.  

One may use these expressions as a basis to 
devise and investigate different, effective stopping criteria in large-scale numerical settings. For example, while $\varepsilon$ appears explicitly 
in Theorem \ref{thm:subspacebound} 
in a way that suggests its use as a relative residual norm threshold; while Theorem \ref{thm:F2r_fro_norm} suggests a
scaling of the residual norm by the norm of the solution vector as another possible 
stopping criterion.
These and related ideas are the focus of on-going work. 
\subsection{Quantities of interest in derived bounds}
\label{sec:condnumber}
By combining Theorem \ref{thm:F2r_fro_norm} with Theorem \ref{thm:Fper_bound}, one observes that 
perturbation effects of the inexact solves on the system theoretical (model reduction related) measures critically depend on the four quantities:  The norm of the oblique projector 
$\widetilde{\bfPhi}_r= \widetilde{\bfV}_r(\widetilde{\bfW}_r^T\widetilde{\bfV}_r)^{-1}\widetilde{\bfW}_r^T$ of the underlying model reduction problem, reciprocals of the minimum singular values of the scaled primitive bases 
$\widetilde{\bfV}_r\bfD$ and 
$\widetilde{\bfW}_r\bfD$; and
 the stopping criterion $\varepsilon$ for the inexact solves, (which affects
 $\displaystyle \max_{i} \frac{\|\bfeta_i\|} {\|\widetilde{\bfv}_i\|}$ and $\displaystyle \max_{i} \frac{\|\bfxi_i\|} {\|\widetilde{\bfw}_i\|}$.)

The $\varepsilon$ term is associated directly with inexact solves and is under the control of the user.  The remaining quantities 
$\varsigma_{\min}(\widetilde{\bfV}_r\bfD)^{-1}$, $\varsigma_{\min}(\widetilde{\bfW}_r\bfD)^{-1}$ and $ \| \widetilde{\bfPhi}_r \|$,
depend largely on the  selection of interpolation points $\{\sigma_i\}$ and tangent directions, but the influence of interpolation data 
on the magnitude of these quantities is difficult to anticipate.

In this section, we will investigate experimentally the effects of
the interpolation point selection on the three quantities of interest, $\varsigma_{\min}(\widetilde{\bfV}_r\bfD)^{-1}$, $\varsigma_{\min}(\widetilde{\bfW}_r\bfD)^{-1}$ and $ \| \widetilde{\bfPhi}_r \|$,  appearing in the derived bounds. 
  These quantities are continuous with respect to the primitive basis vectors, $\{\widetilde{\bfv}_1,\,\cdots,\, \widetilde{\bfv}_r\}$ and
$\{\widetilde{\bfw}_1,\,\cdots,\, \widetilde{\bfw}_r\}$ in neighborhoods where $\widetilde{\bfW}_r^T\widetilde{\bfV}_r$ is nonsingular (i.e., where the projector $\widetilde{\bfPhi}_r$ is well defined). Thus it will be sufficient to examine how the magnitudes of the quantities of interest depend on interpolation data presuming that the necessary linear solves are done \emph{exactly}; for modest convergence thresholds, the effect of inexact solves on these magnitudes is secondary to the effect of interpolation point location.

\begin{figure}[ht]
  \centerline{
\includegraphics[width=12.0cm] {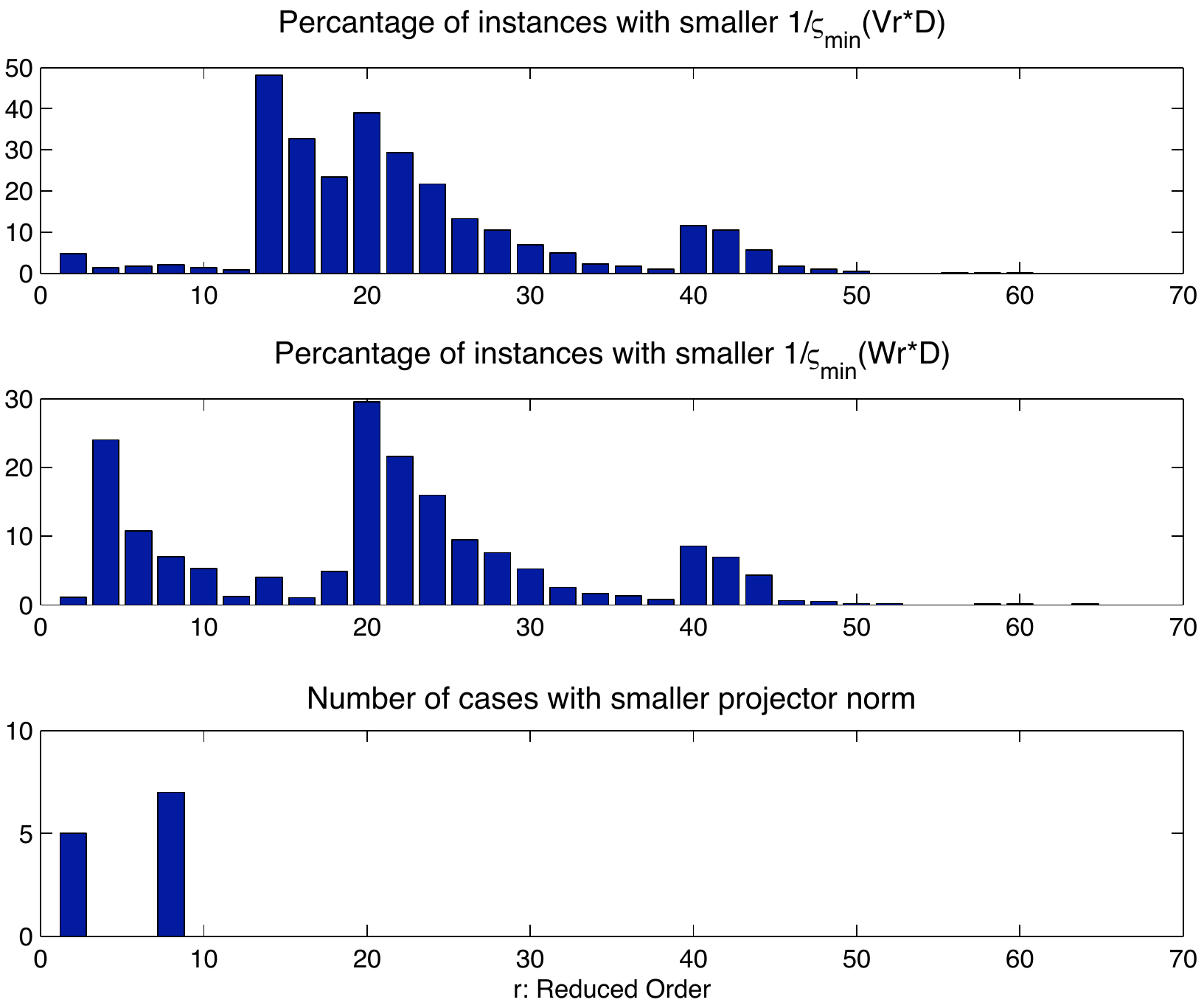}
   }
  \caption{Comparison of  $\varsigma_{\min}({\bfV}_r\bfD)^{-1}$, $\varsigma_{\min}({\bfW}_r\bfD)^{-1}$ and $ \|{\bfPhi}_r \|$
for random shift selections relative to values for $\h2$-optimal shifts }  
  \label{figure:cond}
  \end{figure}

For our numerical study, we use the International Space Station 12A Module as the full-order model. The model has order $n=1412$. We examine a single-input single-output subsystem, $\cbfH(s)$,  reducing the order from $1412$  to order $r$ with $r$ varying from $2$ to $70$ in increments of two. For each reduced order, we chose  $2000$ random shift selections and computed 
 $\varsigma_{\min}({\bfV}_r\bfD)^{-1}$, $\varsigma_{\min}({\bfW}_r\bfD)^{-1}$ and $ \|{\bfPhi}_r \|$.  
For each $r$, $r/2$ shifts were sampled from a uniform distribution on a rectangular region in the positive half-plane: 
$\displaystyle \left\{z\in\IC\left|
\begin{array}{l} \min_\lambda |\mathsf{Re}(\lambda)|\leq\mathsf{Re}(z)\leq \max_\lambda |\mathsf{Re}(\lambda)| \\
|\mathsf{Im}(z)|\leq \max_\lambda |\mathsf{Im}(\lambda)|
\end{array} \right. \right\}$, where the $\max$ and $\min$ are taken over all the poles of the system. The remaining $r/2$ shifts were taken to be the complex conjugates of this random sample, so as to produce a shift configuration that was closed under conjugation.
Additionally for each $r$, we applied model reduction using the $\h2$-optimal interpolation points generated by the method of \cite{gugercin2008hmr}. Then, for each $r$, out of the $2000$ randomly generated shift selections, we counted the number of cases where the random shift selection yielded smaller values of $\varsigma_{\min}({\bfV}_r\bfD)^{-1}$, $\varsigma_{\min}({\bfW}_r\bfD)^{-1}$ and $ \|{\bfPhi}_r \|$. The results are shown in Figure \ref{figure:cond}. Figure \ref{figure:cond}-(a) and -(b) show that for most of the cases, the $\h2$-optimal interpolation points yield smaller values for $\varsigma_{\min}({\bfV}_r\bfD)^{-1}$, $\varsigma_{\min}({\bfW}_r\bfD)^{-1}$. Indeed, for $r\geq 48$,  the $\h2$-optimal points  produced  smaller values in more than $99\%$ of the cases. Also, for the last three cases: $r=66$, $r=68$, and $r=70$, the $\h2$-optimal interpolation points always yielded smaller quantities. The results are even more dramatic for the projector norm, which is important in scaling the perturbation effects caused by inexact solves, see Theorem \ref{thm:F2r_fro_norm}: Out of $70,000$ cases ($2000$ selections for each $r$ value), the $\h2$-optimal interpolation point selection  produced smaller condition numbers in all except $7$ instances:
 $5$ instances for $r=2$, and $2$ instances for $r=8$.  These numerical results illustrate that  $\h2$-optimal interpolation points can be expected to yield smaller  values for $\varsigma_{\min}({\bfV}_r\bfD)^{-1}$, $\varsigma_{\min}({\bfW}_r\bfD)^{-1}$ and $ \|{\bfPhi}_r \|$,  and hence should produce reduced order models that are more robust with respect to perturbations. 

\begin{figure}[ht]
  \centerline{
\includegraphics[width=12.0cm] {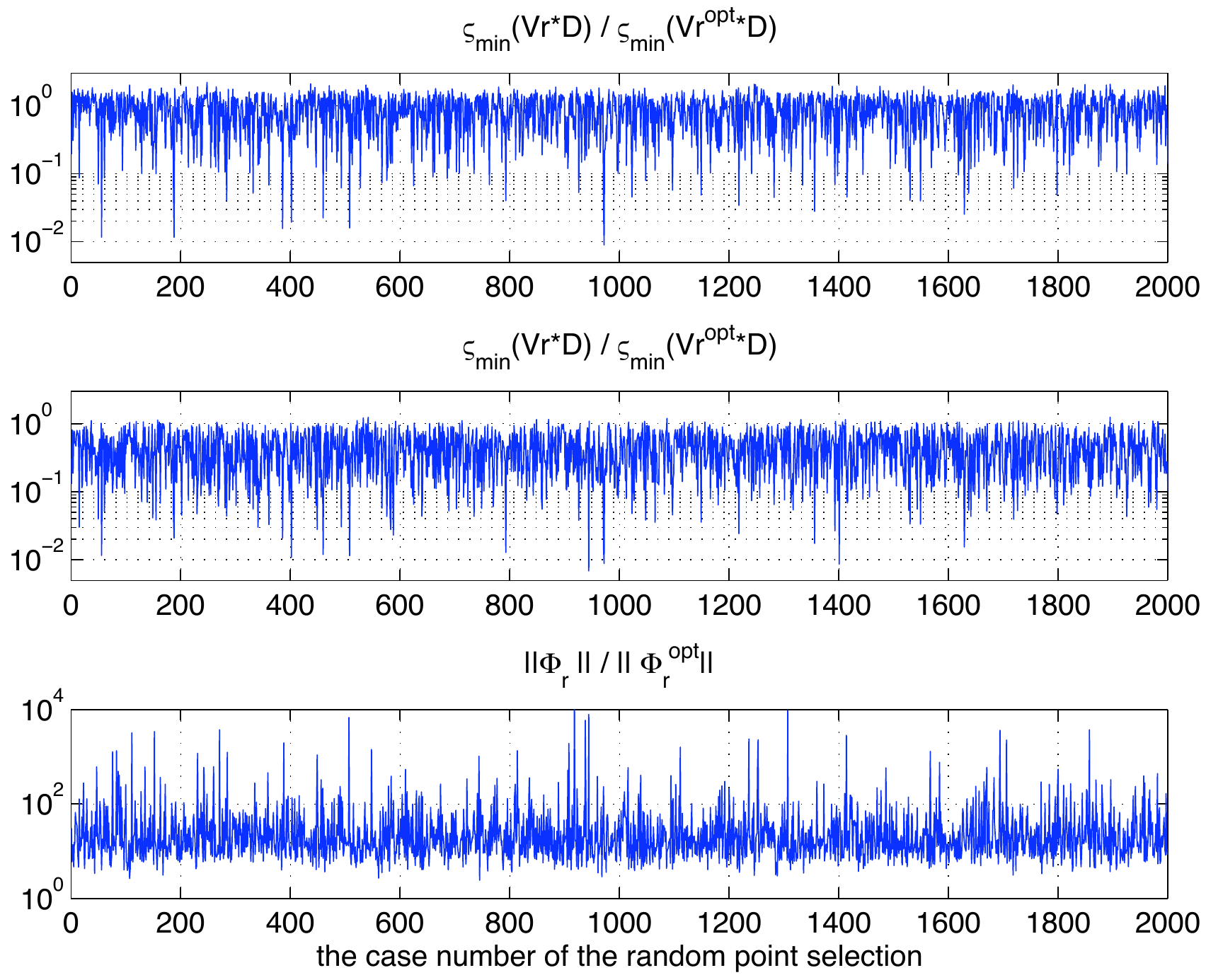}
   }
  \caption{Detailed comparison for $r=14$}  
  \label{fig:smin_compare}
  \end{figure}

Figure \ref{figure:cond} also shows that for $r=14$,  $48\%$ of the randomly selected shifts yielded smaller values of $\varsigma_{\min}({\bfV}_r\bfD)^{-1}$. However, when we inspected the $2000$ randomly selected shift sets for $r=14$ in more detail, we observed some interesting additional features.  We computed the three quantities $\varsigma_{\min}({\bfV}_r\bfD)^{-1}$, $\varsigma_{\min}({\bfW}_r\bfD)^{-1}$ and $ \|{\bfPhi}_r \|$ for each of the $2000$ randomly selected shift sets, and compared them with the corresponding value derived from an $\h2$-optimal shift selection. The results are shown in Figure \ref{fig:smin_compare}. The top plot shows $\varsigma_{\min}({\bfV}_r\bfD) / \varsigma_{\min}({\bfV}_r^{\mbox{\scriptsize\textsf{opt}}}\bfD)$ where ${\bfV}_r^{\mbox{\scriptsize\textsf{opt}}}$ stands for the primitive interpolatory basis for the $\h2$-optimal points. The bigger this ratio, the better the random shift selection. Even though for  $48\%$ of the cases, the random selection was better, the highest this ratio becomes is $2.20$, i.e., the random shifts were never much better than a factor of $2$ better than what $\h2$-optimal shifts provided. 
For the remaining $52\%$ of the cases, the randomly selected shifts were worse, and often worse by a factor of $100$ or more. The situation for $\bfW_r$ is shown in the middle plot.  Once more, the situation is much more drastically in the favor of the $\h2$-optimal interpolation points when the projector norm is inspected; the bottom plot in Figure \ref{fig:smin_compare} which depicts the ratio $\| \bfPhi_r\| / \| \bfPhi_r^{\mbox{\scriptsize\textsf{opt}}}\|$ where $ \bfPhi_r^{\mbox{\scriptsize\textsf{opt}}}$ denotes the projector for the $\h2$-optimal points. As illustrated in Figure \ref{figure:cond}, there are no random shift cases yielding a smaller projector norm. Furthermore, in many cases the projector norm for the random shift selection is almost $4$ order of magnitudes higher than that of the $\h2$-optimal points. Indeed, on average the projector norm for the random points is  $8.19 \times 10^1$ times higher. These numbers change more in the favor of the $\h2$-optimal points as $r$ increases. For example, for $r=50$, while the ratio $\varsigma_{\min}({\bfV}_r\bfD) / \varsigma_{\min}({\bfV}_r^{\mbox{\scriptsize\textsf{opt}}}\bfD)$  becomes only as high as $1.48$, it becomes as low as $2.89 \times 10^{-4}$ for some random selections; Also, the ratio 1 can  reach as high as $2.91\times 10^5$. For $r=70$,  $\| \bfPhi_r\|$ for random selection is $1.73 \times 10^2$ times higher than
$ \| \bfPhi_r^{\mbox{\scriptsize\textsf{opt}}}\|$  on average.    
The three quantities we have been investigating  appear to be extremely well conditioned for $\h2$-optimal interpolation points. Even for $r=70$,  both $\varsigma_{\min}({\bfV}_r^{\mbox{\scriptsize\textsf{opt}}}\bfD)^{-1}$,  $\varsigma_{\min}({\bfW}_r^{\mbox{\scriptsize\textsf{opt}}}\bfD)^{-1}$ 
remain smaller than $10$ and $\| \bfPhi_r^{\mbox{\scriptsize\textsf{opt}}}\|$ is smaller than $7$.


\section{Inexact Solves in Optimal Interpolatory Approximation}  
\label{sec:optimalH2}

The quality of the reduced-order model in interpolatory model reduction clearly depends on the selection of interpolation points and tangent directions. Until recently, this selection process was mostly {\it ad hoc}, and this factor had been the principal  disadvantage of interpolatory model reduction.  For systems in  standard first-order state-space form, 
 Gugercin {\it et  al.} \cite{gugercin2008hmr}
have produced that an $\h2$-optimal interpolation point / tangent direction selection strategy and proposed an Iterative Rational Krylov Algorithm ({\newIRKA}) to generate interpolatory reduced-order models that are (locally) optimal with respect to the $\h2$ norm.  (An $\h2$-optimal interpolation point selection strategy is still unknown for the general coprime factorization framework.)  In this section, we investigate the behavior of inexact solves within the $\h2$-optimal interpolatory approximation setting, specifically examining the behavior when inexact solves are employed in {\newIRKA}. In the rest of this section,  we  briefly review the optimal $\h2$ approximation problem and the method of \cite{gugercin2008hmr}.  We then show how inexact solves can be  employed effectively in this setting and discuss  observed effects on optimality of the  final reduced model.   Our discussion  focuses on
systems in first-order descriptor form:
\begin{equation} \label{desc_form}
\cbfH(s) = \bfC (s \bfE - \bfA)^{-1}\bfB
\end{equation}
where $\bfE,\bfA \in \IR^{n \times n}$, $\bfB \in \IR^{n\times m}$ and 
$\bfC \in \IR^{p\times n}$. 

\subsection{Optimal $\h2$ approximation problem}
Given the full-order system as in (\ref{desc_form}), the goal of the optimal
${\EuScript H}_2$ model reduction problem is to find a reduced-order model $\cbfH_r(s)$ that minimizes the $\h2$ error; i.e. 
\begin{equation} \label{h2opt} 
 \left\| \cbfH - \cbfH_r \right\|_{\h2} =  \min_{\substack{G_r\ \mbox{\scriptsize stable} \\ \dim(G_{r})=r}} \left\| \cbfH - G_{r} \right\|_{\h2}.
\end{equation}

Many researchers have worked on this problem. 
These efforts can be grouped into two categories: Lyapunov-based optimal 
$\h2$ methods
such as 
\cite{yan1999aaa,spanos1992ana,halevi1992fwm,hyland1985top,wilson1970oso,zigic1993contragredient};
and interpolation-based optimal $\h2$ methods such as 
\cite{meieriii1967aol,gugercin2008hmr, gugercin2006rki, vandooren2008hom,bunse-gerstner2009hom, gugercin2005irk, kubalinska2007h0i, beattie2007kbm,beattie2009trm}. 
Here, we will focus on the interpolation-based approach. However we note that 
 Gugercin {\it et al.} \cite{gugercin2008hmr} has shown that these two frameworks are theoretically equivalent; hence motivating the use of  interpolatory approaches to optimal 
$\h2$ approximation since they are numerically superior to the Lyapunov-based approaches.

Since the optimization problem (\ref{h2opt}) is nonconvex, obtaining a global minimizer  is a hard task and can be intractable. The usual approach is to find reduced order models that satisfy first-order necessary optimality conditions. 
Meier and Luenberger \cite{meieriii1967aol}  introduced  
interpolation-based $\h2$-optimality conditions for SISO systems.  Analogous
 $\h2$-optimality conditions for MIMO systems 
 have recently been developed 
by \cite{gugercin2008hmr,bunse-gerstner2009hom,vandooren2008hom} which in turn have led  to analogous algorithms for the MIMO case; 
see \cite{gugercin2008hmr,bunse-gerstner2009hom} for more details. 
\begin{theorem}  \label{h2cond}
Given  $\cbfH(s) =  \bfC (s \bfE - \bfA)^{-1}\bfB$, let $\cbfH_r(s)=\sum_{i=1}^r\frac{1}{s-\widehat{\lambda}_i}\bfc_i\bfb_i^T$
 be  the best $r^{\rm th}$ order 
 approximation of $\cbfH$ with respect to the $\h2$ norm.
Then
\begin{align}
  \mbox{(a)}\ & \cbfH(-\widehat{\lambda}_k) \bfb_k =\cbfH_r(-\widehat{\lambda}_k) \bfb_k,
   \quad \mbox{(b)}\  \bfc_k^T \cbfH(-\widehat{\lambda}_k) =  \bfc_k^T\cbfH_r(-\widehat{\lambda}_k),
   \label{H2optcond}\\
  & \mbox{ and }\quad \mbox{(c)}\  \bfc_k^T \cbfH'(-\widehat{\lambda}_k) \bfb_k 
   =  \bfc_k^T\cbfH_r'(-\widehat{\lambda}_k) \bfb_k 
\quad \mbox{  for } k=1,\,2,\,...,\,r. \nonumber
\end{align}
\end{theorem}

\subsubsection{An algorithm for interpolatory optimal $\h2$ model reduction}
Theorem \ref{h2cond} reveals that any $\h2$ optimal reduced-order model 
$\cbfH_r(s)$ is a bi-tangential Hermite interpolant to $\cbfH(s)$ at mirror images of the reduced-order poles.  However,  since the interpolation points and  the tangent directions (and consequently, $\bfV_r$ and $\bfW_r$), depend on the final reduced-model to be computed, they are not known {\it a priori}. The \emph{Iterative Rational Krylov Algorithm} (\newIRKA) of   \cite{gugercin2008hmr} resolves this problem
by iteratively correcting the interpolation points and the directions as outlined in Algorithm \ref{sucratkry}: The reduced-order order poles are reflected across the imaginary axis to become the next set of interpolation points; the tangent directions are 
corrected using  residue directions from the current reduced model. Upon convergence, the resulting interpolatory reduced-order model 
satisfies the necessary conditions of Theorem  \ref{h2cond}.
For further details on  \newIRKA, see \cite{gugercin2008hmr}.

\begin{center}
    {\small
    \begin{algorithm} \label{sucratkry} 
    {\bf \newIRKA~for MIMO $\h2$ Optimal Tangential Interpolation}        
    \begin{enumerate}
       \item Make an initial $r$-fold shift selection: $\{\sigma_1,\ldots,\sigma_r\}$   and initial tangent directions 
       $\hat{\bfb}_1,\ldots,\hat{\bfb}_r$ and 
       $\hat{\bfc}_1,\ldots,\hat{\bfc}_r$.
       \item  $\bfV_r =  \left[(\sigma_1 \bfE- \bfA)^{-1}\bfB\hat{\bfb}_1~\cdots~(\sigma_r \bfE - \bfA)^{-1}\bfB\hat{\bfb}_r~\right]$\\
       $\bfW_r =  \left[({\sigma_1}\, \bfE - \bfA^T)^{-1}\bfC^T\hat{\bfc}_1 \,\cdots \,({\sigma_r}\, \bfE - \bfA^T)^{-1}\bfC^T\hat{\bfc}_1~\right]$.
       \item while (not converged)
       \begin{enumerate}
       \item $\bfA_r = \bfW_r^T \bfA \bfV_r$, $\bfE_r = \bfW_r^T \bfE \bfV_r$, $\bfB_r = \bfW_r^T\bfB$, and $\bfC_r = \bfC \bfV_r$
       \item Compute $\bfY^*\bfA_r\bfX = diag(\tilde{\lambda}_i)$ and $\bfY^*\bfE_r\bfX = \bfI_r$
        where $\bfY^*$ and $\bfX$ are\\ the left and right eigenvector matrices for $\lambda\bfE_r-\bfA_r$.
              \item $\sigma_i \longleftarrow -{\lambda_i(\bfA_r,\bfE_r)}$ for $i=1,\ldots,r$, $\hat{\bfb}_i^* \longleftarrow \bfe_i^T\bfY^*\bfB_r$ and 
 $\hat{\bfc}_i \longleftarrow \bfC_r \bfX \bfe_i$.
     \item  $\bfV_r =  \left[(\sigma_1 \bfE - \bfA)^{-1}\bfB\hat{\bfb}_1~\cdots~(\sigma_r \bfE - \bfA)^{-1}\bfB\hat{\bfb}_r~\right]$
       \item $\bfW_r =  \left[({\sigma_1}\, \bfE - \bfA^T)^{-1}\bfC^T\hat{\bfc}_1 \,\cdots \,({\sigma_r}\, \bfE - \bfA^T)^{-1}\bfC^T\hat{\bfc}_1~\right]$. 
       \end{enumerate}
       \item $\bfA_r = \bfW_r^T \bfA \bfV_r$, $\bfE_r = \bfW_r^T \bfE \bfV_r$, $\bfB_r = \bfW_r^T \bfB$, $\bfC_r = \bfC\bfV_r$
       \end{enumerate}
       \end{algorithm}
       }
    \end{center}

\subsection{Inexact Iterative Rational Krylov Algorithm (\newinIRKA)} \label{section:iirka} 
For large system order, one may see from Algorithm \ref{sucratkry}, 
that the main cost of \newIRKA~ 
will generally be solving $2r$ large linear systems at each step. If the \newIRKA~iteration converges 
in $k$ steps, a total of $2rk$ linear systems will need to be solved. 
In settings where system dimension reaches into the millions, 
iterative linear system solvers become necessary and inexact linear system
solves must be incorporated into \newIRKA. 
    We refer to the modified
algorithm as the \emph{Inexact Iterative Rational Krylov
Algorithm} (\newinIRKA) and describe it in Algorithm \ref{inexactirka} below. 
We employ the Petrov-Galerkin framework for the inexact solves. In Algorithm  \ref{inexactirka},
the function $\bfF_{\mathsf{PG}} $ in $$[\tildebfv_i,\tildebfw_i] = \bfF_{\mathsf{PG}}\left(\bfA,\bfE,\bfB,\sigma_i,\bfb_i,\bfc_i,\mathbf{v}^{(0)},\mathbf{w}^{(0)},\varepsilon\right)$$  
denotes an inexact solve using a Petrov-Galerkin framework to approximately solve the linear systems $ (\sigma_i \bfE - \bfA) \bfv_i = \bfB \bfb_i$ and  $ (\sigma_i \bfE - \bfA)^T\, \bfw_i =\bfC^T \, \bfc_i$  with initial guesses 
$\mathbf{v}^{(0)}$ and $\mathbf{w}^{(0)}$, respectively,
and a relative residual termination tolerance $\varepsilon$, i.e., at the end,
$$\frac{\|   (\sigma_i \bfE - \bfA)\tildebfv_i- \bfB \bfb_i\|}{\|\bfB \bfb_i\|}\leq \varepsilon ~~~{\rm and}~~~
\frac{\|   (\sigma_i \bfE - \bfA)^T\tildebfw_i- \bfC^T \bfc_i\|}{\|\bfC^T \bfc_i\|}\leq \varepsilon 
$$.

\begin{center}
    {\small
 \begin{algorithm}  \label{inexactirka} 
    {\bf \newinIRKA~for MIMO $\h2$ Optimal Tangential Interpolation}     
    \begin{enumerate}
\item Make an initial $r$-fold shift selection: $\{\sigma_1,\ldots,\sigma_r\}$   and initial tangent directions 
       $\hat{\bfb}_1,\ldots,\hat{\bfb}_r$ and 
       $\hat{\bfc}_1,\ldots,\hat{\bfc}_r$.
\item
 $[\tildebfv_i,\tildebfw_i] = \bfF_{\mathsf{PG}}\left(\bfA,\bfE,\bfB,\sigma_i,\bfb_i,\bfc_i,\mathbf{0},\mathbf{0},\varepsilon\right)$  for $i=1,\ldots,r$
\item $~\widetilde{\bfV}_r = \left[~\tildebfv_1,~\tildebfv_2,~\ldots,~\tildebfv_r~\right]$~~and~~ 
$\widetilde{\bfW}_r = \left[~\tildebfw_1,~\tildebfw_2,~\ldots,~\tildebfw_r~\right]$.
\item while (not converged)
\begin{enumerate}
 \item $\widetilde{\bfA}_r = \widetilde{\bfW}_r^T \bfA \widetilde{\bfV}_r$, $\widetilde{\bfE}_r = \tildebfWr^T \bfE \tildebfVr$, $\widetilde{\bfB}_r = \tildebfWr^T\bfB$, and $\widetilde{\bfC}_r = \bfC \tildebfVr$
  \item Compute $\bfY^*\widetilde{\bfA}_r\bfX = diag(\tilde{\lambda}_i)$ and $\bfY^*\widetilde{\bfE}_r\bfX = \bfI_r$
        where $\bfY^*$ and $\bfX$ are the left and right eigenvector matrices of $\lambda\bfE_r-\bfA_r$.
              \item $\sigma_i \longleftarrow -{\lambda_i(\widetilde{\bfA}_r,\widetilde{\bfE}_r)}$ for $i=1,\ldots,r$, $\hat{\bfb}_i^* \longleftarrow \bfe_i^T\bfY^*\widetilde{\bfB}_r$ and \\
 $\hat{\bfc}_i \longleftarrow \widetilde{\bfC}_r \bfX \bfe_i$.
\item
 $[\tildebfv_i,\tildebfw_i] = \bfF_{\mathsf{PG}}\left(\bfA,\bfE,\bfB,\sigma_i,\bfb_i,\bfc_i,\tildebfv_i,\tildebfw_i,\varepsilon\right)$  for $i=1,\ldots,r$
\item $~\widetilde{\bfV}_r = \left[~\tildebfv_1,~\tildebfv_2,~\ldots,~\tildebfv_r~\right]$ and 
$\widetilde{\bfW}_r = \left[~\tildebfw_1,~\tildebfw_2,~\ldots,~\tildebfw_r~\right]$.
\end{enumerate}
 \item $\widetilde{\bfA}_r = \widetilde{\bfW}_r^T \bfA \widetilde{\bfV}_r$, $\widetilde{\bfE}_r = \tildebfWr^T \cbfE \tildebfVr$, $\widetilde{\bfB}_r = \tildebfWr^T\bfB$, and $\widetilde{\bfC}_r = \bfC \tildebfVr$
\end{enumerate}
\end{algorithm}
     }
    \end{center}
    
As discussed and illustrated in \cite{gugercin2008hmr,antoulas2010imr}, in most cases \newIRKA~converges rapidly; that is, the interpolation 
points and directions at the $k^{\rm th}$ step of \newIRKA~stagnate rapidly with respect to $k$.  Let $\sigma_i^{\rm (k)}$ and $\bfb_i^{(k)}$ denote the $i^{\rm th}$ interpolation point and right-tangential direction, respectively, at the $k^{\rm th}$ step.
Then we expect that as $k$ increases, the solution 
$\bfv_i^{(k)}$ of the linear system 
$(\sigma_i^{(k)}\bfE - \bfA)\bfv^{(k)}=\bfB \bfb_i^{(k}$ from the $k^{\rm 
th}$ step approaches to the solution  $\bfv_i^{(k+1)}$ of the linear system
$(\sigma_i^{(k+1)}\bfE - \bfA)\bfv^{(k+1)}=\bfB \bfb_i^{(k+1)}$at the $(k+1)^{\rm st}$ step. This is precisely the reason that in
Step $4.(d)$ of Algorithm \ref{inexactirka}, we use  $\bfv_i^{(k)}$
  as an initial guess in solving
 $(\sigma_i^{\rm (k+1)}\bfE - \bfA)\bfv^{(k+1)}=\bfB \bfb_i^{(k+1)}$at the $(k+1)^{\rm st}$. We expect that this initialization strategy will speed-up the convergence of the iterative  solves. 
  
  The development of effective stopping criteria based rationally on system theoretic error measures as we have introduced them here is the focus of on-going work.
   Similar approaches toward the design of effective preconditioning techniques and reuse of preconditioners tailored for interpolatory model reduction and especially for optimal $\h2$ approximation are also under  investigation. 

\subsection{Effect of Inexact Solves in the \newinIRKA~Setting}
The first question to answer in \newinIRKA~is whether a statement can be made about the optimality as in the exact \newIRKA~case. Employing the Petrov-Galerkin framework
makes this possible:
\begin{corollary} \label{cor:inirka}
Let $\tildecbfH_r(s)$ be obtained by Algorithm \ref{inexactirka}. Then  $\tildecbfH_r(s)$ satisfies  the necessary conditions for optimal $\h2$ approximation of a near-by full-order model $\tildecbfH(s) = \bfC (s \bfE - (\bfA+\bfF_{2r}))^{-1}\bfB$ where 
$\bfF_{2r}$ is the  rank-$2r$ perturbation matrix defined in (\ref{E2r}).
\end{corollary}
Corollary \ref{cor:inirka} shows that with the help of the underlying Petrov-Galerkin framework, we state that the final reduced model of {\newinIRKA} is an optimal $\h2$ approximation to a nearby full-order model.  

As we discussed in Section \ref{sec:condnumber}, for a good selection
of interpolation points,  interpolatory model reduction is expected to be
robust with respect to perturbations due to inexact solves. Hence, if one feeds
the optimal interpolation points from \newIRKA~into 
an inexact interpolation framework, we expect that the resulting reduced model will be close to the optimal reduced model of \newIRKA. However, the optimal interpolation points are  not known initially and \newinIRKA~will be initiated with 
a nonoptimal initial shift selection.  If the initial interpolation points and directions are poorly selected, at the early stages of the iteration, perturbations due to inexact solves
might be magnified by this poor selection. One can avoid this scenario by using a small termination threshold $\varepsilon$ in the 
early steps of \newinIRKA, and then gradually increase $\varepsilon$ as
the iteration starts to converge. However, we note that
in our numerical experiments using random initialization
strategies,  \newinIRKA~performed robustly and
yielded high fidelity reduced models that are also close to
the true optimal reduced model. This is illustrated in \S \ref{sec:inexactirka} below.  Effective initialization strategies are discussed in \cite{gugercin2008hmr}
as well.

\subsection{Numerical results for \newinIRKA}
\label{sec:inexactirka}
Here we illustrate the usage of inexact solves in the optimal $\h2$ approximation setting 
by  comparing \newIRKA~ with {\newinIRKA}. We use the example of  \S \ref{examples}, but with a finer discretization leading to a state-space dimension of $n=20209$.  We  focus on a MIMO version  using  $2$-inputs and $2$-outputs.

\begin{figure}[ht]
  \centerline{
\includegraphics[width=13.5cm] {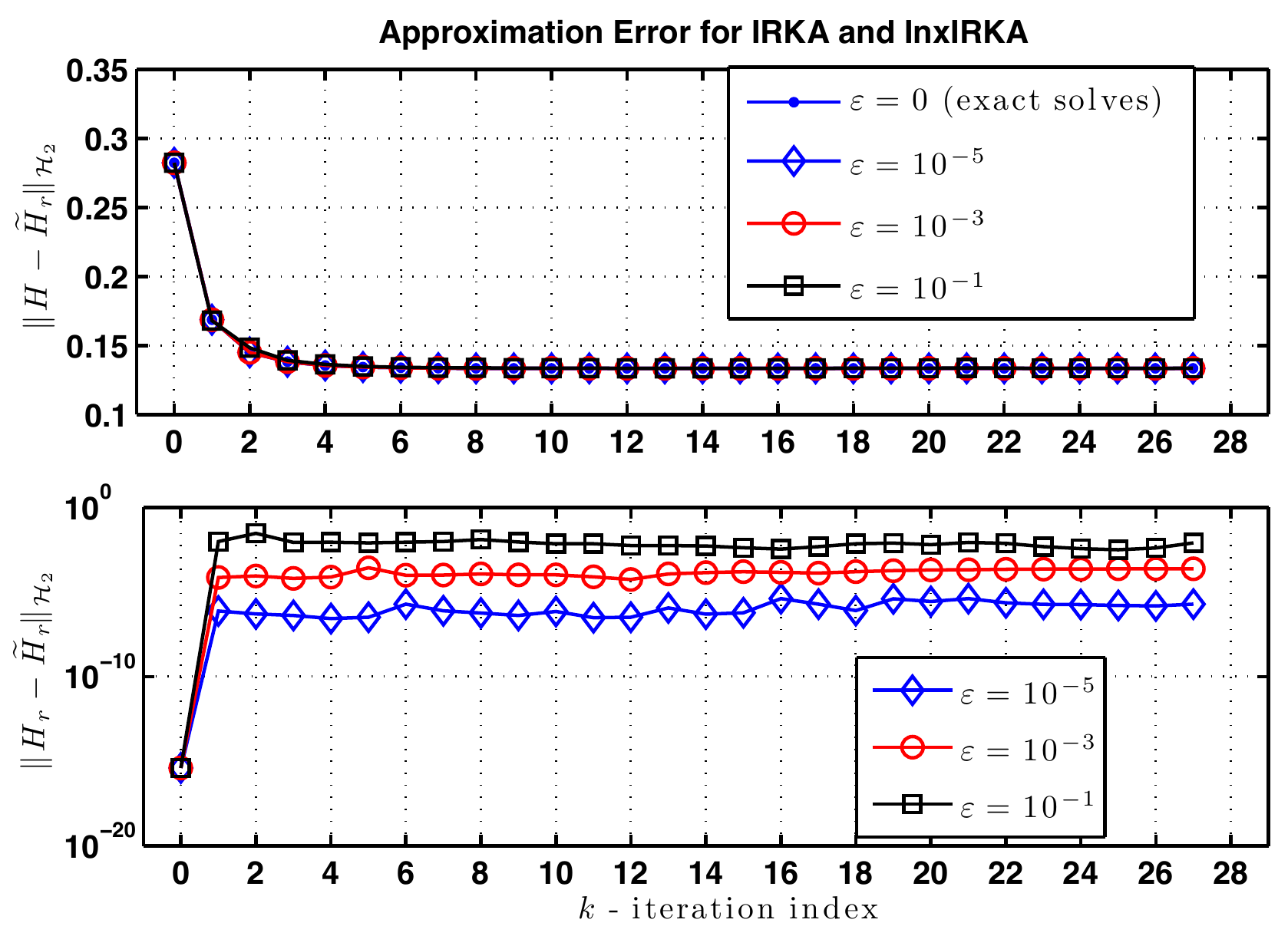}
   }
  \caption{Evolution of the $\h2$ error during {\newIRKA} and {\newinIRKA}}  
  \label{fig:h2evolution}
  \end{figure}

 We reduce the order to $r=6$ using both 
 {\newIRKA} and {\newinIRKA}. In {\newinIRKA}, 
the dual linear systems are
 solved in a Petrov-Galerkin framework using BiCG \cite{barrettTemplates1994} where we use three different values for the relative residual termination threshold of 
 $\varepsilon$:  $ 10^{-5}$, $ 10^{-3}$,  and $ 10^{-1}$.  In all cases, the behavior of {\newinIRKA} is virtually indistinguishable from that of {\newIRKA}.  Starting with the same initial conditions, both  {\newIRKA} and {\newinIRKA} converge within $10$ iteration steps in all $5$ cases. The evolution of the $\h2$
 errors $\| \cbfH - \cbfH_r\|_\h2$ and $ \|\cbfH - \tildecbfH_r\|_\h2$ during the course of {\newIRKA} and {\newinIRKA}, respectively, are depicted in the top plot of Figure \ref{fig:h2evolution}. The figure shows that 
 {\newinIRKA} behavior is almost an exact replica of that of {\newIRKA}. The deviation from the exact 
 {\newIRKA} is noticeable in the graph only for $\varepsilon =  10^{-1}$. To illustrate how much 
 $\cbfH_r $ deviates from $\tildecbfH_r$ as  {\newIRKA} and {\newinIRKA} evolve, we show
 the progress of $\| \cbfH_r - \tildecbfH_r\|_\h2$  in the bottom plot of Figure \ref{fig:h2evolution}. For this example, 
 we initialized both  {\newIRKA} 
 and  {\newinIRKA} with an initial reduced-order model (as opposed to specifying initial interpolation points and tangent directions). Thus, $\cbfH_r = \tildecbfH_r$ initially and no linear solvers are involved in the first ($k=0$) step. One could expect that perturbation errors due to inexact solves might accumulate over the course of the {\newinIRKA} iteration, but this does not appear to be the case as this figure illustrates.  The  magnitude of $\| \cbfH_r - \tildecbfH_r\|_\h2$ 
 remains relatively constant throughout the iteration at a magnitude proportional to the termination criterion.
  
 The resulting 
 $\h2$ and $\hinf$ model reduction errors, $\| \cbfH - \cbfH_r\|_\h2$ and $\| \cbfH - \cbfH_r\|_\hinf$ (with $\cbfH_r$ obtained from {\newIRKA}), versus  $\|\cbfH - \tildecbfH_r\|_\h2$ and $\|\cbfH - \tildecbfH_r\|_\hinf$
  (with $\tildecbfH_r$ obtained from {\newinIRKA}) are given as 
 $\varepsilon$ varies in Table \ref{InIRKA_vs_epsilon} below. 
The row corresponding to $\varepsilon=0$ represents the
 errors due to exact {\newIRKA}.
 \begin{table}[ht]
 \centering
 \begin{tabular}{c|c||c}
 $\varepsilon$ & $\h2$ error & $\hinf$ error  \\ \hline
0                            & $3.708415753 \times 10^{-4}$  & $1.084442854 \times 10^{-2}$ \\
$ 10^{-5}$ & $3.708415754 \times 10^{-4}$  & $1.084425703 \times 10^{-2}$ \\
$ 10^{-4}$ & $3.708415778 \times 10^{-4}$ &  $1.084282001 \times 10^{-2}$ \\
$ 10^{-3}$ & $3.708418102\times 10^{-4}$  &  $1.082437228  \times 10^{-2}$ \\
$ 10^{-2}$ & $3.708621743 \times 10^{-4}$  & $1.064836300  \times 10^{-2}$ \\
$ 10^{-1}$ & $3.716780975 \times 10^{-4}$  & $1.055441476  \times 10^{-2}$
 \end{tabular}
 \caption{Evolution of the model reduction errors as $\varepsilon$ varies}
 \label{InIRKA_vs_epsilon}
 \end{table}
  These numbers demonstrate that employing inexact solves in {\newinIRKA} does not degrade the model reduction performance. We also measure the difference between $\cbfH_r$
  and $\tildecbfH_r$ in both $\h2$ and $\hinf$ norms as $\varepsilon$ varies. These results are tabulated in Table \ref{HrHrtilde_vs_epsilon}:
   \begin{table}[ht]
 \centering
 \begin{tabular}{c|c||c}
  $\varepsilon$ & $ {\| \cbfH_r - \tildecbfH_r\|_\h2} $ &  $ {\| \cbfH_r - \tildecbfH_r\|_\hinf}$  \\ \hline
$ 10^{-5}$  &   $5.1921\times 10^{-9}$   &  $2.7776\times 10^{-7}$  \\
$ 10^{-4}$   &  $5.7156\times 10^{-8}$    & $2.4611\times 10^{-6}$  \\
$ 10^{-3}$   &  $6.3982\times 10^{-7}$   &  $2.1043\times 10^{-5}$   \\
$ 10^{-2}$   &  $5.9277\times 10^{-6}$    & $2.0910\times 10^{-4}$ \\
$ 10^{-1}$   &  $2.2056\times 10^{-5}$   &  $2.9228\times 10^{-3}$  \end{tabular}
 \caption{Evolution of the perturbation error as $\varepsilon$ varies}
 \label{HrHrtilde_vs_epsilon}
 \end{table}
     Note that  while $\| \cbfH - \cbfH_r\|_\h2$ and $\| \cbfH - \cbfH_r\|_\hinf$ are respectively $\mathcal{O}(10^{-4})$ and $\mathcal{O}(10^{-2})$, the contributions attributable to
    $\cbfH_r - \tildecbfH_r$ are much smaller in magnitude and do not alter the resulting (optimal) model reduction performance in any significant way.  If one  were to convert the perturbation errors in Table \ref{HrHrtilde_vs_epsilon}
    to  relative error (as opposed to the displayed absolute error), both
    $ {\| \cbfH_r - \tildecbfH_r\|_\h2} $ and $ {\| \cbfH_r - \tildecbfH_r\|_\hinf}$ starts at $\mathcal{O}(10^{-6})$ 
    for  $\varepsilon = 10^{-5}$, and increases linearly by one order as $\varepsilon$ increases by the same amount.
    
    We finally list, in Table \ref{optsigmas}, the final exact and inexact optimal interpolation points due to 
    {\newIRKA}, and {\newinIRKA} for $\varepsilon = 10^{-3}$ and $\varepsilon = 10^{-1}$:
   \begin{table}[ht]
    $$
    \begin{array}{c|c|c}
    \sigma_i({\newIRKA})~~&~~\sigma_i({\newinIRKA}), \varepsilon =1\times 10^{-3} & 
    ~~\sigma_i({\newinIRKA}), \varepsilon =1\times 10^{-1}
    \\ \hline
    1.0802\times10^{-5}  ~~&~~1.0800\times10^{-5}&~~1.2396\times10^{-5}\\
   9.7164\times 10^{-4} ~~&~~ 9.7080\times10^{-4}&~~ 9.5860\times10^{-4}\\
   6.6310\times 10^{-3}   ~~&~~ 6.6246\times10^{-3}&~~ 6.5923\times10^{-3}\\
   5.7925\times 10^{-2}  ~~&~~5.7938\times10^{-2}&~~5.7929\times10^{-2}\\
   9.0460\times 10^{-1}   ~~&~~9.0419\times10^{-1}&~~8.9877\times10^{-1}\\
   1.4127\times 10^0~~   ~~&~~ 1.4126\times 10^0~~&~~ 1.4104\times 10^0~~
   \end{array}
    $$
    \caption{Optimal interpolations points as $\varepsilon$ varies}
   \label{optsigmas}
    \end{table}
   Not surprizingly, the resulting interpolation points are very close to each other (though not the same).   This can be viewed as another illustration of the fact that $\tildecbfH_r$ is an $\h2$ optimal approximation to a nearby full-order system.

\begin{figure}[ht]
  \centerline{
\includegraphics[width=13.5cm] {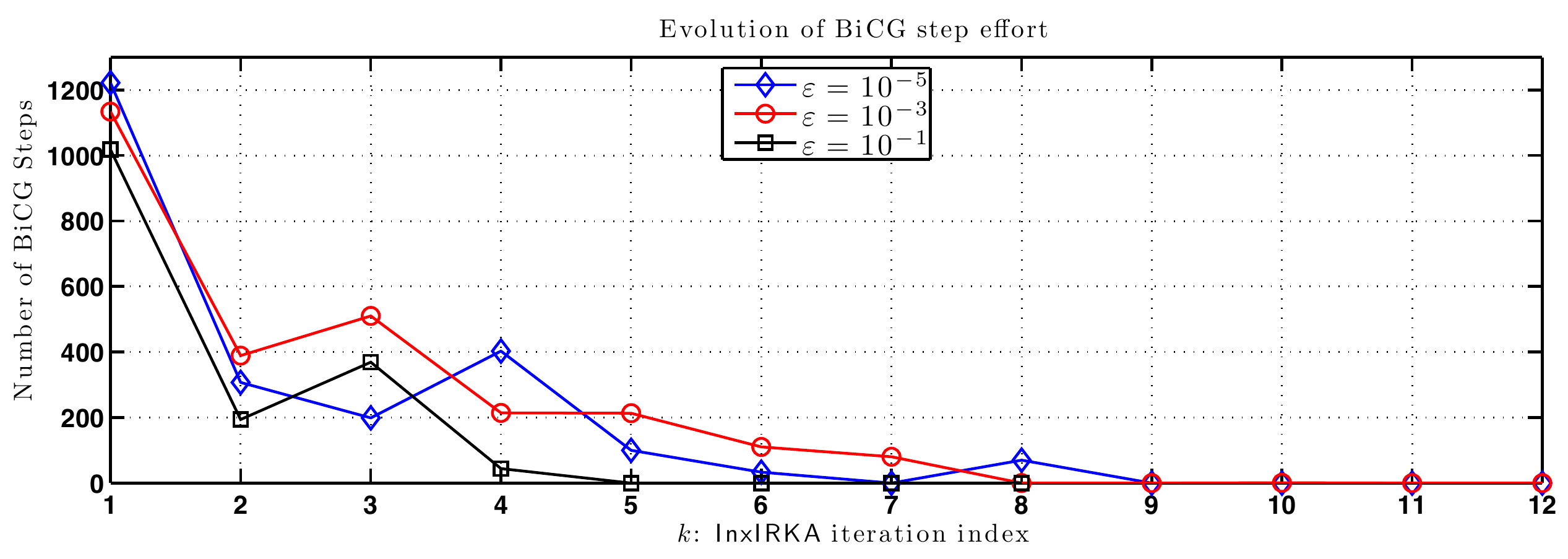}
   }
  \caption{Evolution of BiCG effort during {\newinIRKA} for shift closest to the imaginary axis}  
  \label{fig:bicgcount}
  \end{figure}

As discussed above,  in the implementation of {\newinIRKA}, we  used the solution vectors from the previous step as
the initial guess for the linear system in the next step taking advantage of the convergence in the interpolation points and tangent directions.
To illustrate the effectiveness of this simple approach,
throughout {\newinIRKA} we monitor the number of BiCG steps required to solve
each linear system. We illustrate the behavior only for one of the interpolation points. We choose the interpolation points closest to the  imaginary axis since 
these produce the hardest linear systems to solve and invariably contribute most to the cost of inexact solves. Figure \ref{fig:bicgcount} depicts the the number of BiCG steps required as {\newinIRKA} proceeds for these interpolation points using 
three different stopping criteria  $\varepsilon = 10^{-5}$, $\varepsilon = 10^{-3}$ and $\varepsilon = 10^{-1}$.  The figure clearly illustrates that re-using the solutions from the previous steps works very effectively in reducing the overall cost of the BiCG. The number of BiCG steps goes from $1200$ down to $200$ in $3$ to 
$4$ steps.
  
\section{Structure-preserving interpolation for descriptor systems}
\label{sec:struct}
The backward error analysis of \S \ref{sec:backwarderror} has been presented for the transfer functions in the generalized coprime 
factorization form as in (\ref{HdecompD}). In this section, we show that stronger conclusions on the structure of the reduced system can be drawn in the case the system has a realization as a descriptor system, that is,
\begin{equation} \label{Hfirst}
\cbfH(s) = \bfC (s\bfE -\bfA)^{-1} \bfB 
\end{equation}
where $\bfE,\bfA \in \IR^{n \times n}$, $\bfB \in \IR^{n \times m}$,  and  $\bfC \in \IR^{p \times n}$ are constant matrices.
In this case,  for the interpolation points $\left\{\sigma_j\right\}_{j=1}^r$, and  the tangent directions
 $\left\{\bfb_j\right\}_{j=1}^r$ and $\left\{\bfc_j\right\}_{j=1}^r$, 
 the associated  primitive interpolatory bases 
 $\bfV_r$  and $\bfW_r$  can be obtained from 
 (\ref{k1}) and  (\ref{k2}) using $\cbfK(s) = s \bfE - \bfA$, $\cbfB(s) = \bfB$ (constant matrix) and $\cbfC(s) = \bfC$ (constant matrix). 
Then, the resulting reduced-order model is given by 
\begin{equation} \label{Hredfirst}
\cbfH_r(s) = \bfC_r (s\bfE_r -\bfA_r)^{-1} \bfB_r  
\end{equation}
where 
\begin{equation} \label{ArBrCr}
\bfE_r = \bfW_r^T\bfE \bfV_r,~~
\bfA_r =  \bfW_r^T\bfA \bfV_r,~~
\bfB_r = \bfW_r^T\bfB,~~{\rm and}~~
\bfC_r =  \bfC \bfV_r.
\end{equation}
Let the set $\mathcal{S} = \{\sigma_i,\bfb_i,\bfc_i\}$ denote given tangential interpolation data. 
Define the matrices $\IL[\cbfH,\mathcal{S}] \in \IC^{r \times r}$ and 
$\IM[\cbfH,\mathcal{S}]\in \IC^{r \times r}$ corresponding to the dynamical system $\cbfH(s)$ and  interpolation data  $\mathcal{S}$:
\begin{equation} \label{defineL}
\left(\IL[\cbfH,\mathcal{S}]\right)_{i,j} :=
\left\{  \begin{array}{ll}
\displaystyle 
  \frac{\bfc_i^T\left(\cbfH(\sigma_i)-\cbfH(\sigma_j)\right)\bfb_j}{\sigma_i-\sigma_j}  & {\rm if}~~ i \neq j \\ \\
\bfc_i^T  \cbfH'(\sigma_i)\bfb_i & {\rm if}~~ i = j 
\end{array} \right.
\end{equation}
\begin{equation}\label{defineM}
\left(\IM[\cbfH,\mathcal{S}]\right)_{i,j} := 
\left\{  \begin{array}{ll}
\displaystyle 
  \frac{\bfc_i^T\left(\sigma_i\cbfH(\sigma_i)-\sigma_j\cbfH(\sigma_j)\right)\bfb_j}{\sigma_i-\sigma_j}  &  {\rm if}~~i \neq j \\ \\
\bfc_i^T  \left. [s\cbfH(s)]'\right|_{s=\sigma_i} \bfb_i & {\rm if}~~i= j 
\end{array} \right.
\end{equation}
 $\IL[\cbfH,\mathcal{S}]$ is the \emph{Loewner matrix} associated with the interpolation data $\mathcal S$ and the dynamical system $\cbfH(s)$, $\IM[\cbfH,\mathcal{S}]$ is the \emph{shifted Loewner matrix}
  associated with the interpolation data $\mathcal S$ and the system $s\cbfH(s)$, see \cite{antoulas2010imr, mayo2007fsg}. The next theorem presents
  a canonical structure for the exact interpolatory reduced-order model (\ref{Hredfirst})-(\ref{ArBrCr}).
\begin{theorem} \cite{mayo2007fsg} \label{thm:mayo}
Given a full-order model $\cbfH(s) = \bfC (s\bfE -\bfA)^{-1}\bfB$ and tangential interpolation data $\mathcal{S} = \{\sigma_i,\bfb_i,\bfc_i\}$,
then the reduced-order quantities in (\ref{ArBrCr}) satisfy 
\begin{equation}
\begin{array}{c}
\begin{array}{cc}
\begin{array}{l}
\bfE_r = -\IL[\cbfH,\mathcal{S}],\\[.1in]
 \bfA_r = -\IM[\cbfH,\mathcal{S}],
\end{array} 
&
\bfB_r = \left[ \begin{array}{c} 
\bfc_1^T\cbfH(\sigma_1)\\ \vdots \\ \bfc_r^T\cbfH(\sigma_r) \end{array}
\right], 
\end{array} \\[.4in]
\bfC_r = [~\cbfH(\sigma_1)\bfb_1,~\ldots,~\cbfH(\sigma_r)\bfb_r~].
\end{array}
\end{equation}
\end{theorem}
\subsection{The Petrov-Galerkin framework and structure preservation}
Theorem \ref{thm:mayo} presents a canonical form for the exact bitangential Hermite interpolant in the case of standard state-space model.
Next we  show that if a Petrov-Galerkin framework is employed in the solution of  the linear systems, the inexact reduced-model will have exactly the same form as the exact one. The result is a direct consequence
of Theorems \ref{backwardHrtilde} and \ref{thm:mayo}.
\begin{corollary} \label{cor:backward}
Given the standard full-order model $\cbfH(s) = \bfC (s\bfE -\bfA)^{-1}\bfB$ together with the 
the interpolation data $\mathcal{S} = \{\sigma_i,\bfb_i,\bfc_i\}$, 
 let the inexact solutions 
$\widetilde{\bfv}_j$ for $(\sigma_j\bfE - \bfA)^{-1}\bfB\bfb_j$ and $\widetilde{\bfw}_j$ for 
$(\sigma_j\bfE - \bfA)^{-T}\bfC^T\bfc_j$ be obtained in a Petrov-Galerkin framework as in (\ref{eqn:pgframework}).
Let $\widetilde{\bfV}_r$ and $\widetilde{\bfW}_r$
denote the corresponding inexact Krylov bases as in (\ref{VrWrtildebases}). Define the residuals 
$$
\widetilde{\bfeta}_j  = (\sigma_j \bfE - \bfA) \widetilde{\bfv}_j-\bfB \bfb_j ~~~~~~{\rm and}~~~~~~\widetilde{\bfxi}_j = 
 (\sigma_j \bfE - \bfA)^T \widetilde{\bfw}_j-\bfC^T\bfc_j.
$$
Let the residual matrices $\mathbf{R}_{\mathbf{b}}$ and $\mathbf{R}_{\mathbf{c}}$, and the rank $2r$ matrix $\bfF_{2r}$ be as defined in
(\ref{eqn:RbRc}) and (\ref{E2r}), respectively.
Then, the inexact interpolatory reduced-order model 
\begin{equation} \label{Htilderedfirst}
\tildecbfH_r(s) = \widetilde{\bfC}_r (s\widetilde{\bfE}_r -\widetilde{\bfA}_r)^{-1} \widetilde{\bfB}_r 
\end{equation}
is an exact Hermite bitangential interpolant for the perturbed full-order model
\begin{equation}   \label{Hperturbedfirstorder}
\tildecbfH(s) = \bfC (s\bfE - (\bfA +\bfF_{2r}))^{-1}
\bfB.
\end{equation}
Morever, the reduced-order quantities satisfy
\begin{eqnarray}\label{ArBrCrtilde}
\begin{array}{c}
\widetilde{\bfE}_r = -\IL[\tildecbfH,\mathcal{S}], \\
 \widetilde{\bfA}_r = -\IM[\tildecbfH,\mathcal{S}],
\end{array}~ 
\widetilde{\bfB}_r = \left[ \begin{array}{c} 
\bfc_1^T\tildecbfH(\sigma_1)\\ \vdots \\ \bfc_r^T\tildecbfH(\sigma_r) \end{array}
\right],~
\widetilde{\bfC}_r = [~\tildecbfH(\sigma_1)\bfb_1,~\ldots,~\tildecbfH(\sigma_r)\bfb_r~].
\end{eqnarray}
where $\IL[\tildecbfH,\mathcal{S}]$ and 
$\IM[\tildecbfH,\mathcal{S}]$ are the Loewner matrices associated with the 
dynamical systems $\tildecbfH(s)$ and $s\tildecbfH(s)$ respectively, and the interpolation data  $\mathcal S$ as defined in (\ref{defineL}) and (\ref{defineM}).
\end{corollary}

Corollary \ref{cor:backward} reveals that the inexact reduced-order model quantities have exactly the same structure as their exact counterparts.  The interpolation data $\mathcal S$ is the same in both cases; the only difference is that $\cbfH(s)$ is replaced by $\tildecbfH(s)$ in the construction that yields the Loewner-matrix structure. The preservation of this structure is independent of the accuracy to which the linear systems are solved. 
In the case where $\bfE = \bfI$, the structure of the exact and inexact reduced-models becomes even simpler:
\begin{corollary} \label{prim_inter_struc}
Assume the hypotheses of Theorem \ref{thm:mayo} with $\bfE = \bfI$.
Then the exact interpolant $\cbfH_r(s) = \bfC_r(s\bfI_r - \bfA_r)^{-1}\bfB_r$ satisfies
 \begin{equation}
 \bfA_r = \bfSigma - \bfQ \bfB,~~~\bfB_r = \bfQ,~~~{\rm and}~~~\bfC_r = [~\bfH(\sigma_1)\bfb_1,~\ldots,~\bfH(\sigma_r)\bfb_r~] \end{equation}
 where
  \begin{equation} \label{defQ}
  \bfQ = (\bfW_r^T\bfV_r)^{-1} \bfW_r^T\bfB,~~\bfSigma = {\rm diag}(\sigma_1,\ldots,\sigma_r)~~~{\rm and}~~~\bfB = [\bfb_1,\ldots,\bfb_r].
 \end{equation}
  Assume the hypothses of Corollary \ref{cor:backward} with $\bfE = \bfI$. 
Then, the inexact interpolant 
$\tildecbfH_r(s) = \widetilde{\bfC}_r (s\bfI_r -\widetilde{\bfA}_r)^{-1} \widetilde{\bfB}_r$ 
 satisfies
 \begin{equation}
 \widetilde{\bfA}_r = \bfSigma - \widetilde{\bfQ} \bfB,~~~\widetilde{\bfB}_r = \widetilde{\bfQ},~~~{\rm and}~~~\widetilde{\bfC}_r = [~\tildecbfH(\sigma_1)\bfb_1,~\ldots,~\tildecbfH(\sigma_r)\bfb_r~] 
 \end{equation}
 where
  \begin{equation} \label{defQtilde}
 \widetilde{\bfQ} = (\widetilde{\bfW}_r^T\widetilde{\bfV}_r)^{-1} \widetilde{\bfW}_r^T\bfB,
  \end{equation}
$\tildecbfH(s)$ is the perturbed full-order model as in (\ref{Hperturbedfirstorder})
with $\bfE = \bfI$, and $\bfSigma$ and $\bfB$ are as defined in (\ref{defQ}).
\end{corollary}

\noindent
Corollary \ref{prim_inter_struc} illustrates that in the case of $\bfE = \bfI$,
both of the reduced system matrices, $\bfA_r$ and $\widetilde{\bfA}_r$, are perturbations of rank $\min(r,m,p)$
to the diagonal matrix of interpolation points, $\bfSigma$.

\bibliographystyle{plain}
\bibliography{DOEbiblio}

\begin{thebibliography}{10}

\bibitem{Kap09}
Kapil Ahuja.
\newblock Recycling {B}i-{L}anczos algorithms: {BiCG}, {CGS}, {BiCGSTAB}.
\newblock Master's thesis, Virginia Tech, Blacksburg, Virginia, August 2009.

\bibitem{antoulas2005approximation}
A.C. Antoulas.
\newblock {\em {Approximation of Large-Scale Dynamical Systems (Advances in
  Design and Control)}}.
\newblock Society for Industrial and Applied Mathematics Philadelphia, PA, USA,
  2005.

\bibitem{antoulas2010imr}
A.C. Antoulas, C.A. Beattie, and S.~Gugercin.
\newblock Interpolatory model reduction of large-scale dynamical systems.
\newblock In J.~Mohammadpour and K.~Grigoriadis, editors, {\em Efficient
  Modeling and Control of Large-Scale Systems}. Springer-Verlag, 2010.

\bibitem{barrettTemplates1994}
R.~Barrett, M.~Berry, T.F. Chan, J.~Demmel, J.M. Donato, J.~Dongarra,
  V.~Eijkhout, R.~Pozo, C.~Romine, and H.~Van~der Vorst.
\newblock {\em {Templates for the Solution of Linear Systems: Building Blocks
  for Iterative Methods}}.
\newblock Society for Industrial Mathematics, 1994.

\bibitem{beattie2007kbm}
C.A. Beattie and S.~Gugercin.
\newblock {Krylov}-based minimization for optimal $\mathcal{H}_2$ model
  reduction.
\newblock {\em 46th IEEE Conference on Decision and Control}, pages 4385--4390,
  Dec. 2007.

\bibitem{beattie2009ipm}
C.A. Beattie and S.~Gugercin.
\newblock {Interpolatory projection methods for structure-preserving model
  reduction}.
\newblock {\em Systems \& Control Letters}, 58(3):225--232, 2009.

\bibitem{beattie2009trm}
C.A. Beattie and S.~Gugercin.
\newblock A trust region method for optimal {$\mathcal{H}_2$} model reduction.
\newblock {\em 48th IEEE Conference on Decision and Control}, Dec. 2009.

\bibitem{benner2004slc}
P.~Benner.
\newblock Solving large-scale control problems.
\newblock {\em Control Systems Magazine, IEEE}, 24(1):44--59, 2004.

\bibitem{benner2004ens}
P.~Benner and J.~Saak.
\newblock Efficient numerical solution of the {LQR}-problem for the heat
  equation.
\newblock {\em Proc. Appl. Math. Mech}, 4(1):648--649, 2004.

\bibitem{bunse-gerstner2009hom}
A.~Bunse-Gerstner, D.~Kubalinska, G.~Vossen, and D.~Wilczek.
\newblock {$\mathcal{H}_2$-optimal model reduction for large scale discrete
  dynamical MIMO systems}.
\newblock {\em Journal of Computational and Applied Mathematics}, 2009.
\newblock doi:10.1016/j.cam.2008.12.029.

\bibitem{gallivan2003mrv}
K.~Gallivan, A.~Vandendorpe, and P.~Van~Dooren.
\newblock {Model reduction via truncation: an interpolation point of view}.
\newblock {\em Linear {A}lgebra and {I}ts {A}pplications}, 375:115--134, 2003.

\bibitem{glover1984aoh}
K.~Glover.
\newblock All optimal {Hankel}-norm approximations of linear multivariable
  systems and their {$L^\infty$}-error bounds.
\newblock {\em International Journal of Control}, 39(6):1115--1193, 1984.

\bibitem{grimme1997kpm}
E.~Grimme.
\newblock {\em {Krylov} Projection Methods for Model Reduction}.
\newblock PhD thesis, Coordinated-Science Laboratory, University of Illinois at
  Urbana-Champaign, 1997.

\bibitem{gugercin2005irk}
S.~Gugercin.
\newblock An iterative rational {Krylov} algorithm {(IRKA)} for optimal
  $\mathcal{H}_2$ model reduction.
\newblock In {\em Householder Symposium XVI}, Seven Springs Mountain Resort,
  PA, USA, May 2005.

\bibitem{gugercin2006rki}
S.~Gugercin, A.C. Antoulas, and C.A. Beattie.
\newblock {A rational Krylov iteration for optimal $\mathcal{H}_2$ model
  reduction}.
\newblock In {\em Proceedings of MTNS}, volume 2006, 2006.

\bibitem{gugercin2008hmr}
S.~Gugercin, A.C. Antoulas, and C.A. Beattie.
\newblock $\mathcal{H}_2$ model reduction for large-scale linear dynamical
  systems.
\newblock {\em SIAM Journal on Matrix Analysis and Applications},
  30(2):609--638, 2008.

\bibitem{halevi1992fwm}
Y.~Halevi.
\newblock {Frequency weighted model reduction via optimal projection}.
\newblock {\em Automatic Control, IEEE Transactions on}, 37(10):1537--1542,
  1992.

\bibitem{hyland1985top}
D.~Hyland and D.~Bernstein.
\newblock The optimal projection equations for model reduction and the
  relationships among the methods of {Wilson}, {Skelton}, and {Moore}.
\newblock {\em Automatic Control, IEEE Transactions on}, 30(12):1201--1211,
  1985.

\bibitem{korvink2005obc}
J.G. Korvink and E.B. Rudnyi.
\newblock {Oberwolfach benchmark collection}.
\newblock In {\em Dimension reduction of large-scale systems: proceedings of a
  workshop held in Oberwolfach, Germany, October 19-25, 2003}, page 311.
  Springer Verlag, 2005.

\bibitem{kubalinska2007h0i}
D.~Kubalinska, A.~Bunse-Gerstner, G.~Vossen, and D.~Wilczek.
\newblock {$\mathcal{H}_2$-optimal interpolation based model reduction for
  large-scale systems}.
\newblock In {\em Proceedings of the $16^{\rm th}$ International Conference on
  System Science}, Poland, 2007.

\bibitem{liu1989spa}
Y.~Liu and B.D.O. Anderson.
\newblock Singular perturbation approximation of balanced systems.
\newblock {\em International Journal of Control}, 50(4):1379--1405, 1989.

\bibitem{mayo2007fsg}
A.J. Mayo and A.C. Antoulas.
\newblock {A framework for the solution of the generalized realization
  problem}.
\newblock {\em Linear Algebra and Its Applications}, 425(2-3):634--662, 2007.

\bibitem{meieriii1967aol}
L.~Meier~III and D.~Luenberger.
\newblock Approximation of linear constant systems.
\newblock {\em Automatic Control, IEEE Transactions on}, 12(5):585--588, 1967.

\bibitem{moore1981pca}
B.~Moore.
\newblock Principal component analysis in linear systems: Controllability,
  observability, and model reduction.
\newblock {\em Automatic Control, IEEE Transactions on}, 26(1):17--32, 1981.

\bibitem{mullis1976som}
C.~Mullis and R.~Roberts.
\newblock Synthesis of minimum roundoff noise fixed point digital filters.
\newblock {\em Circuits and Systems, IEEE Transactions on}, 23(9):551--562,
  1976.

\bibitem{spanos1992ana}
J.T. Spanos, M.H. Milman, and D.L. Mingori.
\newblock A new algorithm for {$L_2$} optimal model reduction.
\newblock {\em Automatica (Journal of IFAC)}, 28(5):897--909, 1992.

\bibitem{szyld2006tmp}
D.B. Szyld.
\newblock {The many proofs of an identity on the norm of oblique projections}.
\newblock {\em Numerical {A}lgorithms}, 42(3):309--323, 2006.

\bibitem{vanderSluis1969}
A.~van~der Sluis.
\newblock {Condition numbers and equilibration of matrices}.
\newblock {\em Numerische Mathematik}, 14(1):14--23, 1969.

\bibitem{vandooren2008hom}
P.~van Dooren, K.A. Gallivan, and P.A. Absil.
\newblock {$\mathcal{H}_2$-optimal model reduction of MIMO systems}.
\newblock {\em Applied Mathematics Letters}, 2008.

\bibitem{wilson1970oso}
DA~Wilson.
\newblock Optimum solution of model-reduction problem.
\newblock {\em Proc. IEE}, 117(6):1161--1165, 1970.

\bibitem{yan1999aaa}
W.Y. Yan and J.~Lam.
\newblock An approximate approach to {$\mathcal{H}_2$} optimal model reduction.
\newblock {\em Automatic Control, IEEE Transactions on}, 44(7):1341--1358,
  1999.

\bibitem{zigic1993contragredient}
D.~Zigic, LT~Watson, and C.~Beattie.
\newblock Contragredient transformations applied to the optimal projection
  equations.
\newblock {\em Linear {A}lgebra and {I}ts {A}pplications}, 188:665--676, 1993.

\end{thebibliography}

\end{document}